\documentclass[11pt]{article}
\usepackage{amsmath,amsxtra}
\usepackage{mathtools}
\usepackage{amssymb}
\usepackage{mathrsfs}
\usepackage{amsfonts}
\usepackage{esint}
\usepackage{amscd}
\usepackage[dvipdf]{graphics}
\usepackage[dvips]{graphicx}
\usepackage{graphicx,subfigure}
\usepackage{latexsym}
\usepackage{epstopdf}
\usepackage{latexsym}
\usepackage{amsmath}
\newcommand{\beqa}{\begin{eqnarray}}
\newcommand{\eeqa}{\end{eqnarray}}
\newcommand{\ba}{\begin{eqnarray*}}
\newcommand{\ea}{\end{eqnarray*}}

\date{}
\newtheorem{definition}{Definition}

\def\sym{{\,\rm sym \,}}

\def\be{{\begin{equation}}}
\def\ee{{\end{equation}}}
\def\Om{{\Omega}}
\def\na{{\nabla}}
\def\Ga{{\Gamma}}
\def\pl{{\partial}}

\def\beq{\arraycolsep=1.5pt\begin{eqnarray}}
\def\eeq{\end{eqnarray}}

\newfont{\Blackboard}{msbm10 scaled 1200}

\newfont{\roma}{cmr10 scaled 1200}

\def\<{{\langle}}
\def\>{{\rangle}}

\newtheorem{thm}{{}\hskip\parindent Theorem}[section]
\newtheorem{lem}{{}\hskip\parindent Lemma}[section]
\newtheorem{pro}{{}\hskip\parindent Proposition}[section]

\newtheorem{cor}{{}\hskip\parindent Corollary}[section]
\newtheorem{dfn}{{}\hskip\parindent Definition}[section]
\newtheorem{rem}{{}\hskip\parindent Remark}[section]

\def\be{\begin{equation}}
\def\ee{\end{equation}}
\def\beq{\arraycolsep=1.5pt\begin{eqnarray}}
\def\eeq{\end{eqnarray}}

\large
\title{\bf The linearized Kirchhoff theory for plates with incompatible prestrain}
\date{}

\author{Yizhao Qin$^a$\thanks{Corresponding author.\ Email: qinyz@mail.tsinghua.edu.cn}\quad
Peng-Fei Yao$^b$\quad \\[0.3cm]
$^a$ Department of Mathematical Science \\
Tsinghua University, Beijing 100084, People¡¯s Republic of China
\\
$^b$ Key Laboratory of  Systems and Control\\
Institute of Systems Science,
Academy of Mathematics and Systems Science\\
Chinese Academy of Sciences, Beijing 100190, P. R. China\\
School of Mathematical Sciences\\
University of Chinese Academy of Sciences, Beijing 100049, China
}

\begin{document}
\maketitle
\footnote{This work is  supported by the National
Science Foundation of China, grants  no. 61473126 and no. 61573342, and Key Research Program of Frontier Sciences, CAS, no. QYZDJ-SSW-SYS011.}
\begin{quote}
\begin{small}
{\bf Abstract} \,\,\, In this paper, we derive a linearized Kirchhoff model from three dimensional nonlinear elastic energy of plates with incompatible prestrain as its thickness $h$ tends to zero and its elastic energy scales like $h^{\beta}$ with $2<\beta<4.$ The incompatible prestrain is given as a Riemannian metric $G(x')$ in the three dimensional thin plate which only depends on mid-plate of the thin plates. The problem is studied rigorously by using a variational approach and establishing the $\Ga-$ limit of the non-Euclidean version of the nonlinear elasticity functional when the gauss curvature of the mid-plate $(\Om, g=G_{2\times2})$ is always positive, negative or zero.
\\[3mm]
\\[3mm]
\end{small}
\end{quote}

\section{Introduction}
\setcounter{equation}{0}
\hskip\parindent
In this paper, we are concerned with the problem of derivation of the two-dimensional limit model of the three dimensional nonlinear elastic thin plates with incompatible prestrain as its thickness tends to zero. The motivation behind this study comes from applications for thin objects with internal prestrain such as growing tissues and various manufactured phenomena (for instance, polymer gels, atomically thin graphene layers, and plastically strained sheets). Shape formation driven by internal prestrain is a very active area of research which has been tackled by both experimentalists and mathematician using both analytic and numerical methods, see for instance \cite{ESK,KES,KHBS,JM} and so on.

Early attempts for replacing the three-dimensional model of a thin elastic structure with planar mid-surface at rest, by a two-dimensional model, were based on a priori simplifying assumptions on the deformations and the stresses, or the idea of using the thickness as a small parameter to establish a limit model by asymptotic methods when the prestrain $G=id.$ For early research, we refer to \cite{FRS} and reference therein. However, the models obtained by means of asymptotic formalism still require a justification through rigorous convergence results. The first result in this respect, was obtained by Le Dret and his coauthors in \cite{DR1}. They used the variational point of view and proved the $\Ga-$ convergence of the 3-dimensional elastic energies, whose prestrain $G=id,$ to a nonlinear membrane energy, valid for loads of magnitude of order 1. In fundamental papers \cite{FJM,FJM2}, S. M\"{u}ller and his coauthors established the geometric rigidity estimates which is the key ingredient in later research, and obtained the hierarchy of the limit model from the three-dimensional nonlinear elasticity energy according to the exponents $\beta\geqslant 2$ of scaling of the elastic energy. They rigorously derived the bending model ($\beta=2$), von-K\'{a}rm\'{a}n model ($\beta=4$), the linear model ($\beta>4$) and the novel intermediate models ($2<\beta<4$).

The work in this area has been further extended to limiting theories for thin shells. The first results for thin nonlinear elastic shells were obtained in \cite{DR2} for the scaling $\beta=0.$ The authors obtained a nonlinear membrane shell model. Another study is due to S. M\"{u}ller and his coauthors who analyzed the case $\beta=2.$ This scaling corresponds to a flexural shell model. Further, M. Lewicka and her partners derived the relevant linear theories ($\beta>4$) and the von--K\'{a}rm\'{a}n-like theories ($\beta=4$) in \cite{LMP2}. Subsequently they proceeded to finalize the analysis for elliptic shells in the intermediate range $2<\beta<4$ and obtained the linearized Kirchhoff theory for thin shells in \cite{LMP}. A similar analysis has been performed in the case of developable shells in \cite{HLP} leading to the proof of the collapse of all two-dimensional limiting theories to the linear theory when $\beta>2.$ Following these findings, a conjecture was made in \cite{LP} about the infinite hierarchy of shell models and the various possible limiting scenarios differentiated by rigidity properties of shells. Recently, the intermediate theory for hyperbolic shells in case of $2<\beta<4$ has been established by P. Yao in \cite{Yao}. Yao solved the linear strain equations in non-characteristic region on hyperbolic surfaces and by applying this results he built the asymptotic theories for hyperbolic shells as that in developable shells in \cite{HLP}.

Most recently, there has been a sustained interest in studying similar problems where the shape formation is driven by the internal prestrain caused by growth, swelling, shrinkage or plasticity. The first work to rigorously study non-trivial configurations of thin prestrained flat films was produced in \cite{LP2}. In this paper, Lewicka and her coauthors introduced the variational formulation of this problem for the non-Euclidean version of the nonlinear elastic energy and established its $\Ga-$ limit under the scaling of its elastic energy $\beta=2.$ Moreover, they also obtained a necessary and sufficient conditions for existence of a $W^{2,2}$ isometric immersion of a given 2d metric into $\mathbb{R}^3.$ The prestrain given in their article is of the form:
\begin{gather*}
\begin{pmatrix}
\begin{array}{ll}
(g_{ij})_{2\times2} & \mathbf{0}\\
\mathbf{0} & 1
\end{array}
\end{pmatrix}
\end{gather*}
Later on, they generalized the results in \cite{LP2} to the case that the prestrain is given as a general Riemmanian metric independent of the thickness variable. In \cite{BLS}, they established the $\Ga-$ limit of the nonlinear, non-Euclidean elastic energy $E^h(u)$ with general thickness-independent Riemmanian metric under the scaling of energy $\beta=2$ and obtained the necessary and sufficient condition for $\inf E^h(u)\sim h^2.$ In \cite{L2R}, the $\Ga-$ limit of the nonlinear, non-Euclidean elastic energy $E^h(u)$ with general thickness-independent Riemmanian metric under the scaling of energy $\beta=4$ was obtained and the corresponding scaling analysis of $\inf E^h(u)$ was also performed. Further, building on their previous resuts, they completed the scaling analysis of $E^h(u)$ with a thickness-dependent Riemann metric satisfying certain
general conditions and the derivation of $\Ga-$ limit of the scaled energies $h^{-2n}E^h,$ for all $n\geqslant 1$ in \cite{L2}. In \cite{L1}, the dimension reduction
for oscillatory metrics and the scaling analysis for its elastic energy $E^h(u)$ was obtained. The prescribed incompatibility metric in this paper exhibits a nonlinear dependence on the transversal variable. Moreover, the case that the prestrain metrics are perturbations of the flat metric $id$ was considered in \cite{LMP3, LMP4, LOP}.

In our article, we consider the derivation of the $\Ga-$ limit of thin plates with incompatible prestrain independent of the thickness variable under the assumption that the scaling of its elastic energy $2<\beta<4.$ To construct the recovery sequence, we associate linear strain equations arising in the prestrained elasticity with that on surface without prestain, which are considered in \cite{LMP,HLP,Yao}. The linear strain equations on surface embedded in $\mathbb{R}^3$ has been explored in \cite{LMP,HLP,Yao} and the authors have established the solvability results of linear strain equations on elliptic, developable and hyperbolic surfaces, and the density, matching property of the infinitesimal isometry on these surfaces.
Note that the recovery sequence we construct have a new term $\frac{x_3^2}{2}d_{\epsilon},$ where the vector field $d_{\epsilon}$ is defined as
\begin{equation*}
(\na u_{\epsilon})^Td_{\epsilon}=-(\na\vec{b}_{\epsilon})^T\cdot\vec{b}_{\epsilon},\quad\mbox{and}\quad\<\vec{b}_{\epsilon},d_{\epsilon}\>=0.
\end{equation*}
We make use of this term to eliminate the unexpected terms arising in our deduction and as a result, we need to construct an isometry or a (generalized) mth order infinitesimal isometry with higher regularity.

To overcome this problem in developable or hyperbolic case, we impose higher regularity on $S,$ namely $\mathcal{C}^{2m+2,1}$ for developable surfaces and $\mathcal{C}^{2m+3,1}$ for hyperbolic surfaces, and use the iteration procedure in Theorem 5.2 in \cite{HLP} to produce the (generalized) mth order infinitesimal isometry with the regularity $\mathcal{C}^{3,1}$ other than $\mathcal{C}^{1,1}$ as in \cite{HLP,Yao}.

In elliptic case, we need to derive $\mathcal{C}^{3,\alpha}$ estimates of the solution of the linear strain equations so that we can establish the isometry of class $\mathcal{C}^{3,\alpha}$ on the surface $S.$ To derive the $\mathcal{C}^{3,\alpha}$ estimates of the solution of the linear strain equations, we come up with a more direct approach for it without using the ADN theory \cite{ADN} for elliptic equations. However, we need to generalize the Euclidean version identity
\begin{equation*}
\pl_i\pl_jy^k=\pl_i(\sym\na y)_{kj}+\pl_j(\sym\na y)_{ik}-\pl_k(\sym\na y)_{ij},\quad i,j,k=1,\cdot\cdot\cdot,n
\end{equation*}
to the vector fields on the surface, which presents the second order covariant derivatives of the solutions of linear strain equations by the combination of the components of the solution itself and its first order derivatives and the components of a given 2nd order tensor field and its first order covariant derivatives. We use it to improve the regularity of the solutions of linear strain equations.

Thus, applying these results, we construct the recovery sequence to establish the $\Ga-\sup$ limit, which completes the computation of the $\Ga-$ limit of the 3d nonlinear elastic energy functional.

Now, we start the formulation of our problem and list the main results in our article.

\subsection{Formulation of problem and the main results}
\hskip\parindent
Let $\Omega\subset\mathbb{R}^2$ be an open, bounded, simply connected domain with a $\mathcal{C}^{1,1}$ boundary. For small $h>0,$ we consider thin plates with mid-plate $\Omega,$ given by
\begin{equation*}
\Omega^h=\Omega\times(-\frac{h}{2},\frac{h}{2})=\{x=(x',x_3): x'\in\Omega,|x_3|<\frac{h}{2}\}.
\end{equation*}
Let $G:\bar{\Omega}\rightarrow\mathbb{R}^{3\times 3}$ be a given smooth Riemannian metric in $\Omega^h$ and $G(x',x_3)=G(x'), \forall (x',x_3)\in\Omega^h.$
The scaled elastic energy of the thin prestrained plates is defined by
\begin{equation*}
E^h(u^h)=\frac{1}{h}\int_{\Omega^h}W(\nabla u^hA^{-1})dx,\quad\forall u^h\in W^{1,2}(\Omega^h,\mathbb{R}^3),
\end{equation*}
where $A=\sqrt{G}$ and $W$ is the energy density of the scaled elastic energy.
We have the following assumptions on the energy density $W$:
\begin{itemize}
\item $W$ is $\mathcal{C}^2$ regular in a neighborhood of $SO(3);$
\item$\exists c>0, \forall F\in\mathbb{R}^{3\times 3}, R\in SO(3),$ we have
\be\label{eq1.1a}
W(R)=0, W(RF)=W(F), W(F)\geqslant cdist^2(F,SO(3)).
\ee
\end{itemize}
Under above hypotheses, we intend to compute the $\Ga-$ limit of $\frac{1}{e^h}E^h$ with
\be\label{eq1.1}
0<\lim_{h\rightarrow 0}\frac{e^h}{h^{\beta}}<\infty,\quad\mbox{for some}\quad 2<\beta<4.
\ee
According to \cite{BLS}, we have that $\lim_{h\rightarrow 0}h^{-2}\inf E^h(u^h)=0$ is equivalent to
\begin{equation}\label{eq1.2}
R_{1212}=R_{1213}=R_{1223}\equiv 0,\quad\mbox{in}\quad \Omega^h.
\end{equation}
And the condition \eqref{eq1.2} is also equivalent to the following claim:
there exists an isometric immersion $y_0:\Omega\rightarrow\mathbb{R}^3$ such that
\begin{eqnarray}\label{eq1.3}
&&\left\{\begin{array}{lll}
(\nabla y_0)^{T}\nabla y_0=G_{2\times 2}, \\
sym((\nabla y_0)^{T}\nabla\vec{b}_0)=0,\\
\end{array}\right.
\end{eqnarray}
where the vector field $\vec{b}_0$ is defined in terms of $y_0$ as follows:
\begin{equation}\label{eq1.4}
Q_0^{T}Q_0=G, Q_0e_1=\partial_1y_0, Q_0e_2=\partial_2y_0\quad\mbox{and}\quad Q_0e_3=\vec{b}_0,\quad\mbox{with}\quad \det Q_0>0.
\end{equation}
Moreover, the isometric immersion $y_0$ is smooth, so is the vector field $\vec{b}_0$ and $y_0$ is unique up to overall rigid motions. This is a consequence of the observation that under \eqref{eq1.3}, the second fundamental form of the surface $y_0(\Om)$ is uniquely given in terms of $G.$ The second equation in \eqref{eq1.3} follows from the fact that the kernel of each quadratic form $\mathcal{Q}_{2,A}$ coincides with $so(2).$ For the future use, let's remark that, denoting the inverse matrix $G^{-1}=(G^{ij})_{i,j=1,2,3},$ we have:
\be\label{eq1.5}
\vec{b}_0=-\frac{1}{G^{33}}(G^{13}\pl_1y_0+G^{23}\pl_2y_0)+\frac{1}{\sqrt{G^{33}}}\vec{\mathbf{n}},
\ee
where $\vec{\mathbf{n}}$ is the unit normal vector on $y_0(\Om).$

Without lose of generality, throughout our paper, we further assume that the smooth isometric immersion $y_0:\Om\rightarrow\mathbb{R}^3$ is an isometric embedding from $\Om$ to $\mathbb{R}^3,$ which means that $y_0$ is a smooth isometry (diffeomorphism plus isometric immersion) from $\Om$ to $y_0(\Om).$

Now we list the results of the lower bound of the scaled elastic energy $\frac{1}{e^h}E^h$ under the assumption \eqref{eq1.1}
\begin{thm}\label{thm2.1}
Let $u^h\in W^{1,2}(\Omega^h,\mathbb{R}^3)$ satisfying $E^h(u^h)\leqslant Ce^h.$ Then, there exist $c^h\in\mathbb{R}^3,$ $\bar{R}^h\in SO(3)$ such that
\begin{equation*}
y^h(x',x_3)=(\bar{R}^h)^{T}u^h(x',x_3)-c^h\in W^{1,2}(\Omega^1,\mathbb{R}^3)
\end{equation*}
and $y_0,$ $\vec{b}_0$ are defined as above. Moreover, we have the following properties of $y^h$ and $E^h(u^h):$

$(i)$\,\,\, $y^h\rightarrow y_0$ in $W^{1,2}(\Omega^1,\mathbb{R}^3)$ and $\frac{1}{h}\partial_3y^h\rightarrow\vec{b}_0$ in $L^2(\Omega^1,\mathbb{R}^3).$

$(ii)$\,\,\, The scaled average displacements:
\begin{equation}\label{eq2.4a}
V^h(x')=\frac{h}{\sqrt{e^h}}\int^{\frac{1}{2}}_{-\frac{1}{2}}y^h(x',x_3)-(y_0(x')+hx_3\vec{b}_0(x'))dx_3
\end{equation}
converges (up to a subsequence) in $W^{1,2}(\Omega,\mathbb{R}^3)$ to some $V\in W^{2,2}(\Omega,\mathbb{R}^3)$ with
\be\label{eq2.4b}
\sym((\nabla y_0)^{T}\nabla V)=0.
\ee

$(iii)$\,\,\, the scaled tangential strains:
\begin{equation*}
\frac{h}{\sqrt{e^h}}\sym((\nabla y_0)^{T}\nabla V^h)
\end{equation*}
converge in $L^2(\Omega,\mathbb{R}^{2\times2})$ to some $e\in L^2(\Omega,\mathbb{R}^{2\times2}_{sym})$

$(iv)$\,\,\, Further, defining the quadratic forms $\mathcal{Q}_3$ and $\mathcal{Q}_{2,A}$ by:
\begin{eqnarray}
&&\mathcal{Q}_3(F)=D^2W(id_3)(F,F), \nonumber \\
&&\mathcal{Q}_{2,A}(x',F_{2\times2})=\min\{\mathcal{Q}_3(A(x')^{-1}\tilde{F}A(x')^{-1}): \tilde{F}\in\mathbb{R}^{3\times3}\quad\mbox{with}\quad \tilde{F}_{2\times2}=F_{2\times2}\}, \label{eq2.4c}
\end{eqnarray}
we have:
\begin{equation}\label{eq2.5}
\liminf_{h\rightarrow0}\frac{1}{e^h}E^h(u^h)\geqslant\mathcal{I}_{\beta}(V)=\frac{1}{24}\int_{\Omega}\mathcal{Q}_{2,A}(x',(\nabla y_0)^T\nabla\vec{p}+(\nabla V)^T\nabla\vec{b_0})dx',
\end{equation}
where the vector field $\vec{p}\in W^{1,2}(\Omega,\mathbb{R}^3)$ is uniquely associated with $V$ by:
\begin{equation}\label{eq2.6}
(\nabla y_0)^T\vec{p}=-(\nabla V)^T\nabla\vec{b_0},\quad\mbox{and}\quad \langle\vec{p},\vec{b_0}\rangle=0.
\end{equation}
\end{thm}

To construct the recovery sequence in the proof of the upper bound of the $\Ga-$ limit of $\frac{1}{e^h}E^h,$ we need to associate the non-Euclidean linear strain equations in our article with the linear strain equation on surface as that in \cite{LMP,HLP,Yao}. Thus, in our article, we establish the upper bound of $\Ga-$ limit for the nonlinear elastic energy when the gauss curvature of $y_0(\Om)$ is always positive, zero or negative. For the general case that the gauss curvature of $y_0(\Om)$ changes its sign on the surface, the solvability of the linear strain equations even when the prestrain $G=id$ is still an open problem because in this case we need to deal with mixed type partial differential equations and we don't have suitable tools at hand.

Before we give the results of the upper bound of $\Ga-$ limit for $\frac{1}{e^h}E^h,$ we list the definitions and the related assumptions for the three types of surfaces mentioned above.
\begin{dfn}\label{def1.1}
The surface embedded in $\mathbb{R}^3$ is said to be elliptic if the second fundamental form $\Pi$ is strictly positive (or strictly negative) definite up to the boundary.
\end{dfn}
Hence, without lose of generality,  for elliptic surfaces, we have
\begin{equation*}
\forall x\in \bar{S},\quad \forall X\in T_xS,\quad \frac{1}{C}|X|^2\leqslant\Pi(x)(X,X)\leqslant C|X|^2.
\end{equation*}
From above, we have that the Gaussian curvature of the elliptic surfaces is strictly positive for all $x\in \bar{S}.$

For the case that $S=y_0(\Om)$ is elliptic surfaces, we introduce the follow assumptions following \cite{LMP}:

{ $\bf(H_e)$}:\,\, $\Omega\subset\mathbb{R}^2$ is an open, bounded, simply connected domain of class $\mathcal{C}^{4,\alpha}$ for some $0<\alpha<1$ and the sectional curvature of the two dimensional Riemannian manifold $(\Om,g=G_{2\times2})$ is strictly positive on $\bar{\Om}.$

Thus, $S=y_0(\Om)$ is of class $\mathcal{C}^{\infty}$ up to its boundary with a $\mathcal{C}^{4,\alpha}$ regular boundary $\pl S.$

Now, we move on to the hyperbolic surface. Let $M$ is a surface embedded in $\mathbb{R}^3$ and $S\subset M$ is an open set in $M.$
\begin{dfn}\label{def1.2}
A region $S\subset M$ is said to be hyperbolic if its Gaussian curvature $\kappa<0$ for all $x\in\bar{S}.$
\end{dfn}
Following \cite{Yao,Yao3}, we introduce the definition of the noncharacteristic region.
\begin{dfn}\label{def1.3}
A region $S\subset M$ is said to be noncharacteristic if one of the following two conditions is satisfied:

$(i)$:\,\, Let
\begin{equation*}
S=\{\alpha(t,s): (t,s)\in(0,a)\times(0,b)\},
\end{equation*}
where $\alpha:[0,a]\times[0,b]\rightarrow M$ is an imbedding map which is a family of regular curves with two parameters $t,s$ such that
\begin{eqnarray*}
&&\Pi(\alpha_t(t,s),\alpha_t(t,s))\neq 0,\quad\forall(t,s)\in[0,a]\times[0,b], \\
&&\Pi(\alpha_s(0,s),\alpha_s(0,s))\neq 0,\quad\Pi(\alpha_s(a,s),\alpha_s(a,s))\neq 0,\quad\forall s\in[0,b], \\
&&\Pi(\alpha_t(0,s),\alpha_s(0,s)=\Pi(\alpha_t(a,s),\alpha_s(a,s))=0, \quad\forall s\in[0,b].
\end{eqnarray*}
$(ii)$:\,\, Let $\alpha(\cdot,s)$ be a closed curve with the period $a$ for each $s\in[0,b]$ given. Let
\begin{equation*}
S=\{\alpha(t,s): t\in[0,a),\quad s\in[0,b]\},
\end{equation*}
where $\alpha:[0,a)\times[0,b]\rightarrow M$ is an imbedding map if $\alpha(\cdot,b)$ is a closed curve; $\alpha:[0,a)\times[0,b)\rightarrow M$ is an imbedding map if $\alpha(\cdot,b)$ is one point. Moreover, for each $s\in[0,b],$
\begin{equation*}
\Pi(\alpha_t(t,s),\alpha_t(t,s))\neq 0,\quad\forall t\in[0,a].
\end{equation*}
\end{dfn}
Here, we give the assumptions for the case that surface $S$ is hyperbolic.

{ $\bf(H_h)$}:\,\,$\Omega\subset\mathbb{R}^2$ is an open, bounded, simply connected domain of class $\mathcal{C}^{2m+3,1},$ where $m\geqslant 2,$ and $(\Om,g=G_{2\times2})$ is isometric to a noncharacteristic region of class $\mathcal{C}^{2m+3,1}$ on hyperbolic surface $M$ through the isometry $y_0.$

Thus, we see that $S=y_0(\Om)\subset M$ is a noncharacteristic region of class $\mathcal{C}^{2m+3,1}.$

At last, we deal with the developable surfaces and give the precise setting following \cite{HLP}.
Let $T>0,$ let $\Gamma\in\mathcal{C}^{1,1}([0,T],\mathbb{R}^2)$ be an arclength parametrized curve, and let $s^{\pm}\in\mathcal{C}^{0,1}([0,T])$ be positive functions.
We set
\begin{equation*}
N=(\Gamma')^{\perp}\in\mathcal{C}^{0,1}([0,T],\mathbb{R}^2)\quad\mbox{and}\quad \tilde{\kappa}=\Ga''\cdot N\in L^{\infty}(0,T)
\end{equation*}
the unit normal and the curvature of $\Ga.$ We introduce the bounded domain:
\begin{equation*}
M_{s^{\pm}}=\{(s,t): t\in(0,T), s\in(-s^{-}(t),s^{+}(t))\},
\end{equation*}
the mapping $\Phi:M_{s^{\pm}}\rightarrow\mathbb{R}^2$ given by:
\begin{equation*}
\Phi(s,t)=\Ga(t)+sN(t)
\end{equation*}
and the open line segment:
\begin{equation*}
[\Ga(t)]=\{\Ga(t)+sN(t): s\in(-s^{-}(t),s^{+}(t))\}.
\end{equation*}
According to \cite{HLP}, we assume that:
\begin{equation*}
[\Ga(t_1)]\cap[\Ga(t_2)]=\emptyset \quad\forall\quad t_1\neq t_2\in[0,T].
\end{equation*}
This condition will be needed to define a developable isometry $u,$ which maps segments $[\Ga(t)]$ to segments in $\mathbb{R}^3.$

Let the curve $\gamma$ be a line of curvature of surface $S$ and the moving frame $\mathbf{r}=(\gamma',v,n)^T\in W^{1,2}((0,T),SO(3))$ be the Darboux frame on the surface along $\gamma$ satisfying
\begin{equation*}
\mathbf{r}'=\begin{pmatrix}
0 & \tilde{\kappa} & \tilde{\kappa}_n \\
-\tilde{\kappa} & 0 & 0 \\
-\tilde{\kappa}_n & 0 & 0
\end{pmatrix}\cdot \mathbf{r}
\end{equation*}
with intial value $\mathbf{r}(0)=id,$ where we set $\gamma(t)=\int_0^t\gamma'$ and $\tilde{\kappa}_n\in L^2(0,T)$ is the normal curvature of $\gamma$ on $S.$
Define the mapping:
\begin{equation*}
u:\Phi(M_{s^{\pm}})\rightarrow\mathbb{R}^3,\quad u(\Phi(s,t))=\gamma(t)+sv(t) \quad\forall (s,t)\in M_{s^{\pm}}.
\end{equation*}
Then, according to Lemma 2.1 in \cite{HLP}, we conclude that $u$ is well defined and  \\
$u\in W_{loc}^{2,2}(\Phi(M_{s^{\pm}}),\mathbb{R}^3)$ is an isometric immersion. Hence, the mapping $(s,t)\mapsto\gamma(t)+sv(t)$ is a line of curvature parametrization of $S.$ Moreover, from above setting, we see that the Gaussian curvature of $S$ is always zero. For more details about the isometric immersion $u,$ we refer to Lemma 2.1 and the related remarks in \cite{HLP}. Now, we give the definition of the developable surfaces.
\begin{dfn}\label{def1.4}
A surface $S\subset \mathbb{R}^3$ is said to be developable of class $\mathcal{C}^{k,1}$ if there are $\Ga, N, s^{\pm}, \gamma, v, n, \tilde{\kappa}, \tilde{\kappa}_n, \Phi$ and $u$ as above(also as that in Lemma 2.1 in \cite{HLP}) such that $\Phi(M_{s^{\pm}})=\Om,$ $u\in\mathcal{C}^{k,1}(\bar{\Om},\mathbb{R}^3)$ and $S=u(\Om).$
\end{dfn}
Here, we give the assumptions for the case that surface $S=y_0(\Om)$ is developable.

{ $\bf(H_d)$}:\,\, $\Omega\subset\mathbb{R}^2$ is an open, bounded, simply connected domain of class $\mathcal{C}^{2m+2,1},$ where $m\geqslant 2,$ and $(\Om,g=G_{2\times2})$ is isometric to a developable surface $S$ of class $\mathcal{C}^{2m+2,1}$ through the isometry $y_0,$ whose mean curvature is bounded away from zero:
\begin{equation}\label{eq1.5}
tr_g\Pi(p)>0,\quad\forall\quad p\in \bar{S},
\end{equation}
where $g$ is the induced Riemannian metric  on $S\subset\mathbb{R}^3$ of the Euclidean metric in $\mathbb{R}^3.$

Thus, $S=y_0(\Om)\subset\mathbb{R}^3$ is developable of class $\mathcal{C}^{2m+2,1}.$

Now we give the results of the upper bound of $\Ga-$ limit for $\frac{1}{e^h}E^h.$ First, we consider the elliptic case.
\begin{thm}\label{thm1.2}
Assume that the assumptions $\bf(H_e)$ and \eqref{eq1.1} hold. Then for every $V\in W^{2,2}(\Om,\mathbb{R}^3)$ satisfying $\sym((\nabla y_0)^{T}\nabla V)=0,$ there exists $u^h\in W^{1,2}(\Om^h,\mathbb{R}^3)$ with $E^h(u^h)\leqslant Ce^h$ such that $(i),(ii)$ and $(iii)$ of Theorem \ref{thm2.1}
are satisfied with $R^h=id$ and $c^h=0$ and
\begin{equation}\label{eq1.5a}
\limsup_{h\rightarrow0}\frac{1}{e^h}E^h(u^h)\leqslant I_{\beta}(V).
\end{equation}
\end{thm}
Because the solutions of the linear strain equations will lose regularity if the surface $S$ is developable or hyperbolic, we obtain the following results.
\begin{thm}\label{thm1.3}
Assume the assumptions $\bf(H_h)$ or $\bf(H_d)$ and \eqref{eq1.1} hold. And let
\be\label{eq1.6}
e^h=o(h^{\beta_m}), \quad \beta_m=2+\frac{2}{m}
\ee
Then the results in Theorem \ref{thm1.2} hold.
\end{thm}

Finally, combining Theorem \ref{thm2.1} with Theorem \ref{thm1.2} or Theorem \ref{thm1.3}, it's easy to have the following corollary.
\begin{cor}\label{cor1.1}
Let the assumptions of Theorem \ref{thm1.2} or \ref{thm1.3} hold. If $u^h\in W^{1,2}(\Om^h,\mathbb{R}^3)$ is a minimizing sequence to $\frac{1}{e^h}E^h,$ that is
\begin{equation}\label{eq1.7}
\lim_{h\rightarrow0}(\frac{1}{e^h}E^h(u^h)-\inf\frac{1}{e^h}E^h)=0,
\end{equation}
then the appropriate renormalizations $y^h=(\bar{R}^h)^Tu^h(x',hx_3)-c^h\in W^{1,2}(\Om^1,\mathbb{R}^3)$ obey the convergence statements of Theorem \ref{thm2.1} $(i),(ii), (iii)$ and any limit $V$ minimizes the functional $I_{\beta}$ on the space $\{V\in W^{2,2}(\Om,\mathbb{R}^3): \sym((\nabla y_0)^{T}\nabla V)=0\}.$

Moreover, for any (global) minimizer $V$ of $I_{\beta},$ there exists a minimizing sequence $u^h$ satisfying \eqref{eq1.7} together with the convergence results in Theorem \ref{thm1.2}.
\end{cor}
We'll deal with the proof of Theorem \ref{thm1.2} and \ref{thm1.3} in Section 4 later on.

\section{The lower bound of the scaled elastic energy}
\setcounter{equation}{0}
\hskip\parindent
Following the approach of \cite{L2R}, we obtain similar results as in Lemma 2.3 and Corollary 2.4 in \cite{L2R} from the geometric rigidity estimates\cite{FJM}.
\begin{lem}\label{lem2.1}
Let the assumption (\ref{eq1.2}) hold. Then for any $u^h$ satisfying $\lim_{h\rightarrow 0}\frac{1}{h^2}E^h(u^h)=0,$ there exists matrix field $R^h\in W^{1,2}(\Omega,SO(3))$ such that
\begin{eqnarray}
&&\frac{1}{h}\int_{\Omega^h}|\nabla u^h(x)-R^h(x')(Q_0(x')+x_3B_0(x'))|^2dx\leqslant C(E^h(u^h)+h^4); \label{eq2.1}\\
&&\int_{\Omega}|\nabla R^h(x')|^2dx'\leqslant\frac{C}{h^2}(E^h(u^h)+h^4)\label{eq2.2},
\end{eqnarray}
where the matrix field $B_0(x')$ is defined as $B_0e_1=\partial_1\vec{b}_0, B_0e_2=\partial_2\vec{b}_0, B_0e_3=\vec{d}_0$ and the vector field $\vec{d}_0$ is defined as $\langle Q_0^{T}\vec{d}_0,e_1\rangle=-\langle\partial_1\vec{b}_0,\vec{b}_0\rangle, \langle Q_0^{T}\vec{d}_0,e_2\rangle=-\langle\partial_2\vec{b}_0,\vec{b}_0\rangle, \langle Q_0^{T}\vec{d}_0,e_3\rangle=0.$
\end{lem}
Therefore, we consider $u^h$ with $E^h(u^h)\leqslant Ce^h$ under the assumptions (\ref{eq1.2}) and (\ref{eq1.1}) and by Lemma \ref{lem2.1}, we deduce that
\begin{eqnarray}
&&\frac{1}{h}\int_{\Omega^h}|\nabla u^h(x)-R^h(x')(Q_0(x')+x_3B_0(x'))|^2dx\leqslant Ce^h; \label{eq2.3}\\
&&\int_{\Omega}|\nabla R^h(x')|^2dx'\leqslant C\frac{e^h}{h^2}\label{eq2.4}.
\end{eqnarray}

Now, we start the proof of Theorem \ref{thm2.1}.

{\bf Proof of Theorem \ref{thm2.1}}\,\,\,We split the proof into several steps. In the first four steps, we establish the existence of convergent subsequences for the quantities under consideration and obtain properties of their limits. In the final steps, we prove the lower bounds in \eqref{eq2.5}.

Step 1.\,\,\, To prove the claimed convergence properties for $y^h,$ we first set
\begin{equation*}
\bar{R}^h=\mathbb{P}_{SO(3)}\fint_{\Omega^h}\nabla u^h(x)Q_0(x')^{-1}dx.
\end{equation*}
Note that the projection above is is well defined, because for every $x\in\Omega,$ in view of (\ref{eq2.3}):
\begin{eqnarray*}
&&dist^2(\fint_{\Omega^h}\nabla u^h(x)Q_0(x')^{-1}dx,SO(3))\leqslant|\fint_{\Omega^h}\nabla u^h(x)Q_0(x')^{-1}dx-R^h(x')|^2 \\
&&\leqslant C\fint_{\Omega^h}|\nabla u^h-R^h(Q_0+x_3B_0)|^2+C|R^h(x')-\fint_{\Omega}R^hdx'|^2 \\
&&\leqslant Ce^h+C|R^h(x')-\fint_{\Omega}R^hdx'|^2,
\end{eqnarray*}
so that, taking the average on $\Omega,$ by the Poincar\'{e} inequality and (\ref{eq2.4}), we get
\begin{equation*}
dist^2(\fint_{\Omega^h}\nabla u^hQ_0^{-1}dx,SO(3))\leqslant Ce^h+C\int_{\Omega}|\nabla R^h|^2dx'\leqslant C\frac{e^h}{h^2}.
\end{equation*}
In particular, we observe that
\begin{equation}\label{eq2.7}
|\fint_{\Omega^h}\nabla u^hQ_0^{-1}dx-\bar{R}^h|^2\leqslant C\frac{e^h}{h^2}.
\end{equation}
Moreover, by \eqref{eq2.3}, \eqref{eq2.4} and \eqref{eq2.7},
\begin{eqnarray}\label{eq2.8}
&&\fint_{\Omega}|R^h-\bar{R}^h|^2dx'=\fint_{\Omega^h}|R^h-\bar{R}^h|^2dx \nonumber\\
&&\leqslant C\fint_{\Omega^h}(|R^h-\fint_{\Omega}R^h|^2+|\fint_{\Omega}R^h-\fint_{\Omega^h}\nabla u^hQ_0^{-1}|^2)dx+C\fint_{\Omega^h}|\bar{R}^h-\fint_{\Omega^h}\nabla u^hQ_0^{-1}|^2dx\nonumber\\
&&\leqslant C\fint_{\Omega^h}|\nabla R^h|^2dx+C\fint_{\Omega^h}|\nabla u^h-R^h(Q_0+x_3B_0)|^2dx+C\frac{e^h}{h^2}\leqslant C\frac{e^h}{h^2}.
\end{eqnarray}

Let $c^h\in\mathbb{R}^3$ be such that $\int_{\Omega}V^h=0$ where $V^h$ is defined as in \eqref{eq2.4a}. Denote by $\nabla_hy^h$ the matrix whose columns are given by $\partial_1y^h, \partial_2y^h$ and $\frac{\partial_3y^h}{h},$ so that
\begin{equation}\label{eq2.8a}
\nabla_hy^h(x',x_3)=(\bar{R}^h)^T\nabla u^h(x',hx_3).
\end{equation}
Observe that by \eqref{eq2.3} and \eqref{eq2.8},
\begin{eqnarray*}
&&\int_{\Omega^1}|\nabla_hy^h-Q_0|^2dx\leqslant C\fint_{\Omega^h}|\nabla u^h-\bar{R}^hQ_0|^2dx \\
&&\leqslant C\fint_{\Omega^h}|\nabla u^h-\bar{R}^h(Q_0+x_3B_0)|^2dx+C\fint_{\Omega^h}|x_3R^hB_0|^2dx+C\fint_{\Omega^h}|R^h-\bar{R}^h|^2dx \\
&&\leqslant C\frac{e^h}{h^2}.
\end{eqnarray*}
Therefore, $\nabla_hy^h$ converges in $L^2(\Omega^1)$ to $Q_0.$ Observe that the sequence $\{y^h\}$ is bounded in $W^{1,2}(\Omega^1),$ by the choice of $c^h.$ Passing to a subsequence, if necessary, we get that $y^h$ converges weakly in $W^{1,2}(\Omega^1)$ and so, in fact:
\begin{equation*}
y^h\rightarrow y_0\quad\mbox{in}\quad W^{1,2}(\Omega^1,\mathbb{R}^3)\quad\mbox{and}\quad \frac{1}{h}\partial_3y^h\rightarrow \vec{b}_0\quad\mbox{in}\quad L^2(\Omega^1,\mathbb{R}^3).
\end{equation*}

Step 2.\,\,\, Note that, for every $x'\in\Omega,$
\begin{eqnarray}\label{eq2.9}
&&\nabla V^h(x')=(\frac{h}{\sqrt{e^h}}\fint_{-\frac{1}{2}}^{\frac{1}{2}}(\nabla_hy^h(x',x_3)-Q_0(x'))dx_3)_{3\times2} \nonumber\\
&&=\frac{h}{\sqrt{e^h}}(\fint_{-\frac{1}{2}}^{\frac{1}{2}}\nabla_hy^h-(\bar{R}^h)^TR^h(Q_0(x')+hx_3B_0)dx_3)_{3\times2}+\frac{h}{\sqrt{e^h}}(((\bar{R}^h)^TR^h-id)Q_0)_{3\times2}
\nonumber\\
&&\triangleq I_1^h+I_2^h.
\end{eqnarray}
The first term above converges to 0. Indeed,
\begin{eqnarray}\label{eq2.10}
&&\|I_1^h\|^2_{L^2(\Omega)}\leqslant C\frac{h^2}{e^h}\fint_{\Omega^1}|\nabla_hy^h-(\bar{R}^h)^TR^h(Q_0+hx_3B_0)|^2dx \nonumber\\
&&\leqslant C\frac{h^2}{e^h}\fint_{\Omega^h}|\nabla u^h-R^h(Q_0+hx_3B_0)|^2dx\leqslant Ch^2.
\end{eqnarray}
Towards estimating the second term in \eqref{eq2.9}, denote
\begin{equation*}
S^h=\frac{h}{\sqrt{e^h}}[(\bar{R}^h)^TR^h-id].
\end{equation*}
By \eqref{eq2.4} and \eqref{eq2.8}, it follows that
\begin{equation*}
\|S^h\|^2_{L^2(\Omega)}\leqslant C\frac{h^2}{e^h}\int_{\Omega}|R^h-(\bar{R}^h)^T|^2\leqslant C\quad\mbox{and}\quad \|\nabla S^h\|^2_{L^2(\Omega)}\leqslant C\frac{h^2}{e^h}\int_{\Omega}|\nabla R^h(x')|^2dx'\leqslant C.
\end{equation*}
Passing to a subsequence, we can assume that $S^h\rightarrow S$ weakly in $W^{1,2}(\Omega).$ By compact embedding theorem for Sobolev space, we have
\begin{equation}\label{eq2.10a}
S^h\rightarrow S \quad\mbox{in}\quad L^2(\Omega)\quad\mbox{and}\quad S^h\rightarrow S \quad\mbox{in}\quad L^4(\Omega),
\end{equation}
which implies that $I_2^h\rightarrow(SQ_0)_{3\times2}$ in $L^2(\Omega,\mathbb{R}^{3\times2}).$ Consequently, by \eqref{eq2.9}, $\nabla V^h\rightarrow(SQ_0)_{3\times2}$ in $L^2(\Omega,\mathbb{R}^{3\times2}).$ As before, we conclude that $V^h$ converges in $W^{1,2}(\Omega)$ and that its limit $V$ belongs to $W^{2,2}(\Omega,\mathbb{R}^3),$ since
$\nabla V=(SQ_0)_{3\times2}\in W^{1,2}(\Omega).$ We now prove \eqref{eq2.4b}. By definition of $S^h,$
\be\label{eq2.11}
\sym S^h=-\frac{1}{2}\frac{\sqrt{e^h}}{h}(S^h)^TS^h,
\ee
so in view of the boundedness of $\{S^h\}$ in $W^{1,2},$
\begin{equation*}
\|\sym S^h\|_{L^2(\Omega)}\leqslant C\frac{\sqrt{e^h}}{h}\|S^h\|^2_{L^4(\Omega)}\leqslant C\frac{\sqrt{e^h}}{h}\|S^h\|^2_{H^1(\Omega)}\leqslant C\frac{\sqrt{e^h}}{h}.
\end{equation*}
Consequently, $S$ is a skew symmetric field. But $(\nabla y_0)^T\nabla V=(Q_0^TSQ_0)_{2\times2},$ hence \eqref{eq2.4b} follows.

For further use, let's define $\vec{p}\in W^{1,2}(\Omega,\mathbb{R}^3)$ by
\be\label{eq2.12}
[\nabla V|\vec{p}]=SQ_0.
\ee
Since $Q_0^T[\nabla V|\vec{p}]\in so(3),$ it's easily checked that $\vec{p}$ is given solely in terms of $V$ by
\begin{eqnarray}\label{eq2.13}
&&\left\{\begin{array}{lll}
(\nabla y_0)^{T}\vec{p}=-(\nabla V)^{T}\vec{b}_0, \\
\langle\vec{p},\vec{b}_0\rangle=0,\\
\end{array}\right.
\end{eqnarray}

Step 3.\,\,\, We now want to establish convergence in $(iii).$ In view of \eqref{eq2.9} we write
\be\label{eq2.14}
\frac{h}{\sqrt{e^h}}\sym(Q_0^T\nabla V^h)_{2\times2}(x')=\frac{h}{\sqrt{e^h}}\sym(Q_0^TI_1^h)_{2\times2}+\frac{h}{\sqrt{e^h}}\sym(Q_0^TS^hQ_0)_{2\times2}\triangleq J_1^h+J_2^h.
\ee
We first deal with the sequence $J_2^h.$ By \eqref{eq2.10a} and \eqref{eq2.11},
\be\label{eq2.15}
\frac{h}{\sqrt{e^h}}\sym S^h\rightarrow-\frac{1}{2}S^TS=\frac{1}{2}S^2\quad\mbox{in}\quad L^2(\Omega).
\ee
Therefore,
\be\label{eq2.16}
J_2^h\rightarrow-\frac{1}{2}(Q_0^TS^TSQ_0)_{2\times2}=-\frac{1}{2}(\nabla V)^T\nabla V\quad\mbox{in}\quad L^2(\Omega).
\ee
We now turn to $J_1^h.$ Recall that by \eqref{eq2.14}, \eqref{eq2.9} and \eqref{eq2.8a},
\begin{eqnarray}\label{eq2.17}
&&J_1^h=\frac{h}{\sqrt{e^h}}\sym(Q_0^TI_1^h)_{2\times2}=\frac{h^2}{\sqrt{e^h}}\sym(\frac{1}{\sqrt{e^h}}Q_0^T(\bar{R}^h)^T(\fint_{-\frac{1}{2}}^{\frac{1}{2}}\nabla u^h-R^h(Q_0+hx_3B_0)dx_3)_{3\times2}) \nonumber \\
&&=\frac{h^2}{\sqrt{e^h}}\sym(Q_0^T(\bar{R}^h)^T\fint_{-\frac{1}{2}}^{\frac{1}{2}}Z^h(x',x_3)dx_3),
\end{eqnarray}
where the rescaled strains $Z^h$ are defined by
\begin{equation*}
Z^h(x',x_3)=\frac{1}{\sqrt{e^h}}(\nabla u^h(x',hx_3)-R^h(x')(Q_0(x')+hx_3B_0(x')))
\end{equation*}
By \eqref{eq2.3}, $\{Z^h\}$ is bounded in $L^2(\Omega^1,\mathbb{R}^3).$ Thus, up to a subsequence,
\begin{equation*}
Z^h(x',x_3)\rightharpoonup Z\quad\mbox{in}\quad L^2(\Omega^1,\mathbb{R}^3).
\end{equation*}
Thus, we have
\begin{equation*}
J_1^h\rightarrow 0 \quad\mbox{in}\quad L^2(\Omega).
\end{equation*}
The above convergence yields $(iii)$ by \eqref{eq2.14} and \eqref{eq2.16}.

Step 4.\,\,\, We now aim at giving the structure of the weak limit $Z(x',x_3)$ of $Z^h(x',x_3).$ From above, we see that $e=-\frac{1}{2}(\nabla V)^T\nabla V.$
As a tool, we consider the difference quotients $f^{s,h}:$
\begin{equation*}
f^{s,h}(x',x_3)=\frac{1}{s\sqrt{e^h}}(y^h(x',x_3+s)-y^h(x',x_3)-hs(\bar{b}_0+h(x_3+\frac{s}{2})\vec{d}_0))
\end{equation*}
We'll show that $f^{s,h}\rightharpoonup\vec{p}$ in $W^{1,2}(\Omega^1,\mathbb{R}^3),$ as $h\rightarrow0,$ for any $s.$ Write
\begin{equation*}
f^{s,h}(x',x_3)=\frac{1}{\sqrt{e^h}}\fint_0^s\partial_3y^h(x',x_3+t)-h(\bar{b}_0+h(x_3+t)\vec{d}_0)dt,
\end{equation*}
and observe that
\begin{eqnarray*}
&&\frac{1}{\sqrt{e^h}}[\partial_3y^h-h(\bar{b}_0+hx_3\vec{d}_0)]=\frac{h}{\sqrt{e^h}}[\frac{1}{h}\partial_3y^h-(\bar{b}_0+hx_3\vec{d}_0)] \\
&&=\frac{h}{\sqrt{e^h}}[(\bar{R}^h)^T\nabla u^h(x',hx_3)-(Q_0+hx_3B_0)]e_3 \\
&&=\frac{h}{\sqrt{e^h}}(\bar{R}^h)^T[\nabla u^h(x',hx_3)-R^h(Q_0+hx_3B_0)]e_3+S^h(Q_0+hx_3B_0)e_3 \\
&&=h(\bar{R}^h)^TZ^h(x',x_3)e_3+S^h(Q_0+hx_3B_0)e_3.
\end{eqnarray*}
The first term on the right-hand side above converges to 0 in $L^2(\Omega^1),$ because $\{Z^h\}$ is bounded in $L^2(\Omega^1,\mathbb{R}^3),$ while the second term converges to $SQ_0e_3=S\vec{b}_0$ in $L^2(\Omega^1)$ by \eqref{eq2.10a}. Note that $SQ_0e_3=\vec{p}$ by \eqref{eq2.12}. Thus, $f^{s,h}\rightarrow\vec{p}$ in $L^2(\Omega^1).$

We now deal with the derivatives of the studied sequence. Firstly,
\begin{eqnarray*}
&&\pl_3f^{s,h}(x',x_3)=\frac{1}{s\sqrt{e^h}}[\pl_3y^h(x',x_3+s)-\pl_3y^h(x',x_3)-h^2s\vec{d}_0] \\
&&=\frac{1}{s}[\frac{1}{\sqrt{e^h}}(\pl_3y^h(x',x_3+s)-h(\vec{b}_0+h(x_3+s)\vec{d}_0))-\frac{1}{\sqrt{e^h}}(\pl_3y^h(x',x_3)-h(\vec{b}_0+hx_3\vec{d}_0))],
\end{eqnarray*}
converges to 0 in $L^2(\Omega^1,\mathbb{R}^3).$ For $i=1,2,$ the in-plane derivatives read as
\begin{eqnarray*}
&&\pl_if^{s,h}(x',x_3)=\frac{1}{s\sqrt{e^h}}[\pl_iy^h(x',x_3+s)-\pl_iy^h(x',x_3)-hs(\pl_i\vec{b}_0+h(x_3+\frac{s}{2})\pl_i\vec{d}_0)] \\
&&=\frac{1}{s}[(\bar{R}^h)^TZ^h(x',x_3+s)-(\bar{R}^h)^TZ^h(x',x_3)]e_i \\
&&+\frac{1}{s\sqrt{e^h}}[(\bar{R}^h)^TR^h(x')(Q_0(x')+h(x_3+s)B_0(x')) \\
&&-(\bar{R}^h)^TR^h(x')(Q_0(x')+hx_3B_0(x'))]e_i-\frac{h}{\sqrt{e^h}}(B_0e_i+h(x_3+\frac{s}{2})\pl_i\vec{d}_0)\\
&&=\frac{1}{s}[(\bar{R}^h)^TZ^h(x',x_3+s)-(\bar{R}^h)^TZ^h(x',x_3)]e_i+S^hB_0(x')e_i-\frac{h^2}{\sqrt{e^h}}(x_3+\frac{s}{2})\pl_i\vec{d}_0,
\end{eqnarray*}
Hence, by the weak convergence of $Z^h,$
\begin{equation*}
\pl_if^{s,h}(x',x_3)\rightharpoonup \frac{1}{s}[\bar{R}^TZ(x',x_3+s)-\bar{R}^TZ(x',x_3)]e_i+S\vec{b}_0e_i\quad\mbox{in}\quad L^2(\Omega^1,\mathbb{R}^3),
\end{equation*}
where $\bar{R}\in SO(3)$ is an accumulation point of the rotation $\bar{R}^h.$

Consequently, $f^{s,h}\rightharpoonup\vec{p}$ in $W^{1,2}(\Omega^1,\mathbb{R}^3)$ and for $i=1,2:$
\be\label{eq2.18}
s\pl_i\vec{p}=[\bar{R}^TZ(x',x_3+s)-\bar{R}^TZ(x',x_3)]e_i+sS\vec{b}_0e_i,
\ee
which proves that $Z(x',.)e_i$ has linear form and that
\be\label{eq2.19}
(\bar{R}^TZ(x',x_3))_{3\times2}=(\bar{R}^TZ(x',0))_{3\times2}+x_3(\nabla\vec{p}-(SB_0)_{3\times2}).
\ee

Step 5.\,\,\, We now prove the lower bound in $(iv).$ Recall that by the definition of $Z^h,$
\begin{equation*}
\nabla u^h(x',hx_3)=R^h(x')(Q_0(x')+hx_3B_0(x'))+\sqrt{e^h}Z^h(x',x_3).
\end{equation*}
Since $Q_0A^{-1}\in SO(3)$ we have
\begin{equation*}
W(\nabla u^hA^{-1})=W((Q_0A^{-1})^T(R^h)^T\nabla u^hA^{-1})=W(id+h\mathcal{J}+\sqrt{e^h}\mathcal{G}^h),
\end{equation*}
where $\mathcal{J}(x',x_3)=x_3A^{-1}Q_0^TB_0A^{-1}\in so(3)$ and $\mathcal{G}^h=A^{-1}Q_0^T(R^h)^TZ^h(x',x_3)A^{-1}.$
Note that by the weak convergence of $Z^h,$
\begin{equation*}
\mathcal{G}^h\rightharpoonup\mathcal{G}=A^{-1}Q_0^T(\bar{R}^T)Z(x',x_3)A^{-1}\quad\mbox{in}\quad L^2(\Omega^1,\mathbb{R}^{3\times3}).
\end{equation*}
Define the "good set":
\begin{equation*}
\Omega_h=\{x\in\Om^1: (e^h)^{\frac{1}{4}}|G^h|<1\}.
\end{equation*}
By the above, the characteristic function $\chi_{\Om_h}\rightarrow 1$ in $L^1(\Om^1).$ Further, by frame invariance and Taylor expanding of $W$ on $\Om_h:$
\begin{eqnarray*}
&&W(id+h\mathcal{J}+\sqrt{e^h}\mathcal{G}^h)=W(e^{-h\mathcal{J}}(id+h\mathcal{J}+\sqrt{e^h}\mathcal{G}^h)) \\
&&=W(id+\sqrt{e^h}\mathcal{G}^h-\frac{1}{2}h^2\mathcal{J}^2+o(\sqrt{e^h})) \\
&&=\frac{1}{2}\mathcal{Q}_3(\sqrt{e^h}\mathcal{G}^h-\frac{1}{2}h^2\mathcal{J}^2)+o(e^h)+o(\sqrt{e^h})|\sqrt{e^h}\mathcal{G}^h-\frac{1}{2}h^2\mathcal{J}^2| \\
&&+e^h|\mathcal{G}^h
-\frac{h^2}{2\sqrt{e^h}}\mathcal{J}^2+o(1)|^2o(1),
\end{eqnarray*}
where the symbol $o(\cdot)$ refers to the quantities uniformly converge to 0, as $h$ tends to 0. Thus, by H\"{o}lder inequality and dominant convergence theorem,
\beq\label{eq2.20}
&&\liminf_{h\rightarrow0}\frac{1}{e^h}E^h(u^h)\geqslant \liminf_{h\rightarrow0}\frac{1}{e^h}\int_{\Om^1}\chi_{\Om_h}W(id+h\mathcal{J}+\sqrt{e^h}\mathcal{G}^h)dx \nonumber \\
&&\geqslant \liminf_{h\rightarrow0}\frac{1}{2}\int_{\Om^1}\mathcal{Q}_3(\chi_{\Om_h}\sym(\mathcal{G}^h-\frac{h^2}{2\sqrt{e^h}}\mathcal{J}^2))dx+\int_{\Om^1}\chi_{\Om_h}|G^h
-\frac{h^2}{2\sqrt{e^h}}\mathcal{J}^2+o(1)|^2o(1)dx \nonumber \\
&&\geqslant \frac{1}{2}\int_{\Om^1}\mathcal{Q}_3(\sym\mathcal{G})dx.
\eeq
By \eqref{eq2.19}, we have
\beq\label{eq2.21}
&&\mathcal{Q}_3(\sym\mathcal{G})=\mathcal{Q}_3(A^{-1}\sym(Q_0^T(\bar{R}^T)Z(x',x_3))A^{-1}) \nonumber \\
&&=\mathcal{Q}_3(A^{-1}(\sym(Q_0^T(\bar{R}^T)Z(x',0))+x_3\sym(Q_0^T(\nabla\vec{p}-(SB_0)_{3\times2})))A^{-1})
\eeq
Submitting \eqref{eq2.21} into \eqref{eq2.20}, we obtain
\beq\label{eq2.22}
&&\liminf_{h\rightarrow0}\frac{1}{e^h}E^h(u^h)\geqslant\frac{1}{2}\int_{\Om^1}\mathcal{Q}_3(\sym\mathcal{G})dx \nonumber \\
&&=\frac{1}{2}\int_{\Om^1}\mathcal{Q}_3(A^{-1}(\sym(Q_0^T(\bar{R}^T)Z(x',0))A^{-1})dx \nonumber \\
&&+\frac{1}{2}\int_{\Om^1}\mathcal{Q}_3(x_3A^{-1}\sym(Q_0^T(\nabla\vec{p}-(SB_0)_{3\times2}))A^{-1})dx \nonumber \\
&&\geqslant \frac{1}{24}\int_{\Om}\mathcal{Q}_{2,A}(x',\sym((\nabla y_0)^T\nabla\vec{p}+(\nabla V)^T\nabla\vec{b}_0))dx'.
\eeq
Thus, we establish the lower bound. \hfill$\Box$

\section{The linear strain equations and the density, matching properties of infinitesimal isometry}
\setcounter{equation}{0}
\hskip\parindent
To construct the recovery sequence for the upper bound of $\Gamma-$ limit, we need to solve the linear strain equations
\begin{equation}\label{eq3.1}
\sym((\na y_0)^T\na V)=U,
\end{equation}
where $U$ is a second order symmetric tensor field in $\Om.$
By the results of the solvability of the linear strain equations, we have the density and matching properties of the infinitesimal isometry.

Now, we associate the infinitesimal isometry $V$ satisfying \eqref{eq2.4b} in $\Om$ with that on surface $S=y_0(\Om).$ Note that we assume the isometric immersion $y_0$ is an isometry from the two dimensional Riemannian manifold $(\Om,g=G_{2\times2})$ to the surface $S=y_0(\Om)$ embedded in $\mathbb{R}^3.$ Therefore, there exists a the $1-1$ corresponding between the infinitesimal isometry $\tilde{V}$ on $S$ and $V$ in \eqref{eq2.4b} given by the change of variables $V=\tilde{V}\circ y_0.$ Denote $V(x')=\tilde{V}(y_0(x'))$
for $x'\in\Om.$ Let $e_1,e_2\in\mathbb{R}^2$ be the natural base in $\mathbb{R}^2$ and $\tilde{e}_i=\pl_iy_0\in T_{y_0(x')}(S).$
We have
\begin{equation}\label{eq3.2}
(\na y_0)^T\na V(e_i,e_j)=\<\pl_i V,\pl_j y_0\>=\<\na_{\tilde{e}_i} \tilde{V},\tilde{e}_j\>,\quad i,j=1,2,
\end{equation}
and by Gauss formula for submanifolds $\na_{\tilde{e}_i}\tilde{e}_j=D_{\tilde{e}_i}\tilde{e}_j-\Pi(\tilde{e}_i,\tilde{e}_j)\vec{\mathbf{n}},$ where $D$ is the Levi-Civita connection on surface $S$ induced by the Euclidean metric in $\mathbb{R}^3,$ $\Pi$ is the second fundamental form of $S$ and the vector field $\vec{\mathbf{n}}$ represents the unit normal vector on surface $S$ embedding in $\mathbb{R}^3,$
\beq\label{eq3.3}
&&\<\na_{\tilde{e}_i} \tilde{V},\tilde{e}_j\>=\tilde{e}_i\<\tilde{V},\tilde{e}_j\>-\<\tilde{V},\na_{\tilde{e}_i}\tilde{e}_j\> \nonumber \\
&&=\tilde{e}_i\<\tilde{V}^{\top},\tilde{e}_j\>-\<\tilde{V}^{\top},D_{\tilde{e}_i}\tilde{e}_j\>+\<\tilde{V},\vec{\mathbf{n}}\>\Pi(\tilde{e}_i,\tilde{e}_j) \nonumber \\
&&=\<D_{\tilde{e}_i}\tilde{V}^{\top},\tilde{e}_j\>+\<\tilde{V},\vec{\mathbf{n}}\>\Pi(\tilde{e}_i,\tilde{e}_j) \nonumber \\
&&=y_0^{\ast}(D\tilde{V}^{\top})(e_i,e_j)+\<\tilde{V},\vec{\mathbf{n}}\>\cdot y_0^{\ast}\Pi(e_i,e_j),
\eeq
where $\tilde{V}^{\top}$ stands for the projection of $\tilde{V}$ into the tangent bundle $TS.$

Combine \eqref{eq3.2} and \eqref{eq3.3}, we obtain
\be\label{eq3.4}
\sym((\na y_0)^T\na V)=y_0^{\ast}(\sym D\tilde{V}^{\top}+\<\tilde{V},\vec{\mathbf{n}}\>\Pi)=y_0^{\ast}(\sym\na\tilde{V}).
\ee

Denote $U(x')=y_0^{\ast}(\tilde{U}(y_0(x'))).$ Then by \eqref{eq3.4}, we have the $1-1$ corresponding between $V$ satisfying \eqref{eq3.1} and $\tilde{V}$ satisfying
\be\label{eq3.5}
\sym\na\tilde{V}=\tilde{U}\quad\mbox{on}\quad S.
\ee
Thus, we apply the consequences in \cite{LMP, HLP, Yao} of solvability of \eqref{eq3.5}, density and matching properties of infinitesimal isometry on $S$ here.
As we have known, until now, for the solvability of \eqref{eq3.5} on surface $S,$ only the results when the Gaussian curvature of the surface $S$ is always positive, negative or zero are obtained
which means the surfaces are elliptic, hyperbolic or developable, repectively. The general case is still an open problem because in general case we may meet a mixed type partial differential equations on $S,$ for which we don't have effective methods at hand.

Now we cope with solvability of the linear strain equations \eqref{eq3.5}, density and matching properties of infinitesimal isometry on elliptic, developable and hyperbolic surfaces in the following subsections.

\subsection{The results on elliptic surfaces}
\hskip\parindent
For the case that $S=y_0(\Om)$ is an elliptic surface, we arrive at the following results in the frame work in \cite{LMP} and \cite{Yao}.

To cope with the linear strain equations on elliptic surface, we first assume $\Om$ is a simply connected, compact domain with nonempty $\mathcal{C}^{2}$ boundary. From here on in this subsection, $y_0$ is only a parametrization of the surface $S.$ And here we assume $y_0\in\mathcal{C}^{2,1}(\bar{\Om},\mathbb{R}^3)$ and $S$ is parametrized by the single chart $y_0(\Om).$ Denote $\mathcal{V}$ is the infinitesimal isometry space on $S,$ which means that for every $\tilde{V}\in\mathcal{V},$ $\pl_{\tau}\tilde{V}(x)=\tilde{A}(x)\tau, \tilde{A}^T(x)=-\tilde{A}(x)$ for any $x\in S$ and $\tau\in T_xS.$

First, we transfer the linear strain equations to an equivalent system via the methods in \cite{Yao} and we employ the notations there. Let $y=\mathcal{T}\tilde{U}$ be the solution of the linear strain equations $\sym\nabla(\mathcal{T}\tilde{U})=\tilde{U},$ where $\tilde{U}\in L^2(S,T_{sym}^2(S)).$ Let $\{E_1,E_2\}$ be a positively oriented frame field on $S.$
Define a function $v$ and a vector field $u$ as follows:
\beq
&&v=p(y)=\frac{1}{2}[\na y(E_2,E_1)-\na y(E_1,E_2)],\label{eq3.6}\\
&&u=\na y(\vec{\mathbf{n}},E_1)E_1+\na y(\vec{\mathbf{n}},E_2)E_2, \label{eq3.7}
\eeq
where the map $p: W^{1,2}(S,\mathbb{R}^3)\rightarrow L^2(S)$ is a linear operator. And it's easy to see that the definitions of $v,u$ are independent of the choice of a positively oriented orthonormal basis, which means that they are globally defined vector fields on $S.$
By Theorem 2.1 in \cite{Yao}, when the Gaussian curvature $\kappa(p)\neq0$ of $S$ for all $p\in\bar{S}$ we have that the linear strain equations
\begin{equation*}
\sym\nabla(\mathcal{T}\tilde{U})=\tilde{U}
\end{equation*}
are equivalent with the following system satisfied by $u$ and $v:$
\beq
&&\<D^2v,Q^{\ast}\Pi\>=P(\tilde{U})-v\kappa tr_g\Pi+X(v),\quad\mbox{on}\quad S,\label{eq3.8} \\
&&u=Q(\na\vec{\mathbf{n}})^{-1}Q[\Lambda(\tilde{U})-D(tr_g\tilde{U})]-Q(\na\vec{\mathbf{n}})^{-1}Dv, \label{eq3.9}
\eeq
where $D$ is the Levi-Civita connection on $S$ with respect to the induced metric $g=G_{2\times2}$ on $S$ from $\mathbb{R}^3,$ the vector field $\Lambda(\tilde{U})$ is defined by $\<\Lambda(\tilde{U}),\alpha\>=tr_gi_{\alpha}D\tilde{U}$ for $\alpha\in T_pS,$ $p\in S$ and
\begin{eqnarray*}
&&P(\tilde{U})=\<D\{Q[\Lambda(\tilde{U})-D(tr_g\tilde{U})]\},Q^{\ast}\Pi\>-\<Q[\Lambda(\tilde{U})-D(tr_g\tilde{U})],(\na\vec{\mathbf{n}})^{-1}D\kappa\> \\
&&-\kappa tr_g\tilde{U}(Q\na\vec{\mathbf{n}}\cdot,\cdot), \\
&&X=(\na\vec{\mathbf{n}})^{-1}D\kappa.
\end{eqnarray*}
The mapping $Q: T_pS\rightarrow T_pS$ in \eqref{eq3.9} is defined by
\begin{equation*}
Q\alpha=\<\alpha,e_2\>e_1-\<\alpha,e_1\>e_2,\quad\mbox{for all}\quad \alpha\in T_pS.
\end{equation*}
Obviously, $Q\alpha$ is well defined and $Q$ can be generalized to the map $Q: T(S)\rightarrow T(S)$ by $(QX)(p)=QX(p),$ for $p\in S,$ $X\in T(S).$ Moreover, the operator $Q$ further induces an operator on $k(k\geqslant2)$ order tensor fields space $T^k(S),$ denoted by $Q^{\ast}:T^k(S)\rightarrow T^k(S),$ by
\begin{equation*}
(Q^{\ast}T)(X_1,\cdot\cdot\cdot,X_k)=T(QX_1,\cdot\cdot\cdot,QX_k),\quad X_1,\cdot\cdot\cdot,X_k\in T(S),\quad T\in T^k(S).
\end{equation*}
By the definition \ref{def1.1}, we regard the second fundamental form $\Pi$ on $S$ as another Riemannian metric on $S,$ denoted by $\hat{g}.$ From \cite{Yao2}, we have
\be\label{eq3.10}
\<D^2v,Q^{\ast}\Pi\>=\kappa\triangle_{\hat{g}}v+\frac{1}{2\kappa}\Pi(QD\kappa,QDv),
\ee
where $\triangle_{\hat{g}}$ is the Laplacian of the metric $\hat{g}.$ Thus, in elliptic case, equation \eqref{eq3.8} becomes
\be\label{eq3.11}
\triangle_{\hat{g}}v=\frac{1}{\kappa}P(\tilde{U})-vtr_g\Pi+\frac{1}{2\kappa}\tilde{X}(v),
\ee
where $\tilde{X}(v)=2\<(\na\vec{\mathbf{n}})^{-1}D\kappa,Dv\>-\frac{1}{\kappa}\Pi(QD\kappa,QDv).$
Moreover, following \cite{Yao}, we have a more direct relation between $y$ and $u,v:$
\be\label{eq3.12}
\na_{\alpha}y=i_{\alpha}\tilde{U}-vQ\alpha+\<u,\alpha\>\vec{\mathbf{n}},\quad \alpha\in T_pS,\quad p\in S,
\ee
where the inner product for tensor field $i$ is defined by a $(k-1)$ order tensor field $i_XT(X_1,\cdot\cdot\cdot,X_{k-1})=T(X,X_1,\cdot\cdot\cdot,X_{k-1}),$ for $X,X_1,\cdot\cdot\cdot,X_{k-1}\in T(S),$ $T\in T^k(S).$

Now we employ the approaches used in \cite{LMP} and obtain the following results.
\begin{pro}\label{pro3.1}
There exists a linear operator
\begin{equation*}
\mathcal{S}:W^{2,2}(S,T^2(S))\rightarrow W^{2,2}(S)
\end{equation*}
such that
\begin{equation*}
\|\mathcal{S}(\tilde{U})\|_{L^2(S)}\leqslant C\|\tilde{U}\|_{L^2(S)}\quad\mbox{and}\quad\|\mathcal{S}(\tilde{U})\|_{W^{1,2}(S)}\leqslant C\|\tilde{U}\|_{W^{1,2}(S)},
\end{equation*}
and that $v=\mathcal{S}(\tilde{U})$ is a solution to \eqref{eq3.11}, for each $\tilde{U}\in W^{2,2}(S,T^2(S)).$
\end{pro}

By Proposition \ref{pro3.1}, we can obtain the solvability of the linear strain equations.
\begin{thm}\label{thm3.1}
There exists a linear operator
\begin{equation*}
\mathcal{T}: L^2_{sym}(S,\mathbb{R}^{2\times2})\rightarrow\{w\in L^2(S,\mathbb{R}^3): w^{\top}\in W^{1,2}(S)\}
\end{equation*}
such that $\sym\nabla(\mathcal{T}\tilde{U})=\tilde{U},$ for every $\tilde{U}\in L^2(S,T_{sym}^2(S))$ and that
\begin{equation*}
\|(\mathcal{T}\tilde{U})^{\top}\|_{W^{1,2}(S)}+\|\<\mathcal{T}\tilde{U},\vec{\mathbf{n}}\>\|_{L^2(S)}\leqslant C\|\tilde{U}\|_{L^2(S)}
\end{equation*}
\end{thm}

Now we consider the matching property for infinitesimal isometry on $S.$ Because we need the isometry obtained in the following with higher regularity than that in \cite{LMP} for the construction of the recovery sequence, we work in the framework in \cite{LMP} to derive a theorem providing an isometry with higher regularity. First, we derive the $\mathcal{C}^{3,\alpha}$ estimates of $y=\mathcal{T}(\sym((\na\phi)^T\na \psi)).$ Following the methods in \cite{LMP} may also arrive at this results. But we derive this estimates
by a more direct methods without using the ADN thoery(refer to \cite{ADN}).
\begin{pro}\label{pro3.2}
Under the assumptions in Theorem \ref{thm3.2}, one has the following uniform estimates:
\be\label{eq3.13}
\|\mathcal{T}(\sym((\na\phi)^T\na \psi))\|_{\mathcal{C}^{3,\alpha}(\bar{S})}\leqslant C\|\phi\|_{\mathcal{C}^{3,\alpha}(\bar{S})}\|\psi\|_{\mathcal{C}^{3,\alpha}(\bar{S})},
\ee
for all $\phi,\psi\in\mathcal{C}^{3,\alpha}(\bar{S},\mathbb{R}^3).$
\end{pro}
Note that under the assumptions in Theorem \ref{thm3.2}, the surface $S$ is of class $\mathcal{C}^{5,\alpha}$ regularity up to the boundary with $\mathcal{C}^{4,\alpha}$ boundary $\pl S$ and its parametrization $y_0\in\mathcal{C}^{5,\alpha}(\bar{\Om},\mathbb{R}^3).$

Before we tackle the proof of Proposition \ref{pro3.2}, we generalize the Euclidean version identity
\begin{eqnarray*}
&&\pl_i\pl_jy^k=\pl_i(\sym\na y)_{kj}+\pl_j(\sym\na y)_{ik}-\pl_k(\sym\na y)_{ij},\quad i,j,k=1,\cdot\cdot\cdot,n\quad \\
&&\mbox{for}\quad y\in\mathcal{C}^2(\mathbb{R}^n,\mathbb{R}^n)
\end{eqnarray*}
to the vector fields defined on surface $S.$
\begin{lem}\label{lem3.1a}
Let $f\in\mathcal{C}^{2}(\bar{S},\mathbb{R}^3).$ We have the following identity:
\beq
&&\na^2f(E_i,E_j,E_k)=\na(\sym\na f)(E_i,E_j,E_k)+\na(\sym\na f)(E_k,E_i,E_j) \nonumber \\
&&-\na(\sym\na f)(E_j,E_k,E_i)+\<f,\na_{E_j}(D_{E_k}E_i)\>-\<f,\na_{E_i}(D_{E_k}E_j)\> \nonumber \\
&&-\frac{1}{2}[E_j,E_k](\<f,E_i\>)-\frac{1}{2}[E_k,E_i](\<f,E_j\>)-\sym\na f([E_k,E_i],E_j) \nonumber \\
&&+\sym\na f([E_j,E_k],E_i)-\sym\na f([E_i,E_j],E_k)+\frac{1}{2}E_k\<f,[E_i,E_j]\> \nonumber \\
&&+\frac{1}{2}E_j\<f,[E_i,E_k]\>-\frac{1}{2}E_i\<f,[E_j,E_k]\>-E_j(\Pi(E_k,E_i))\<f,\vec{\mathbf{n}}\>-\Pi(E_k,E_i)\Pi(f^{\top},E_j) \nonumber \\
&&+E_i(\Pi(E_j,E_k))\<f,\vec{\mathbf{n}}\>+\Pi(E_k,E_j)\Pi(f^{\top},E_i),\label{eq3.14}
\eeq
where $i,j,k=1,2.$
\end{lem}
{\bf Proof}\,\,\, First we calculate $\na^2f(E_i,E_j,E_k)$ as follows and we set $\na_{\vec{\mathbf{n}}} f=0$ and $f_i=\<f,E_i\>.$
\beq
&&\na^2f(E_i,E_j,E_k)=\na_{E_k}(\na f)(E_i,E_j)=E_k(\na f(E_i,E_j))-\na f(\na_{E_k}E_i,E_j)-\na f(E_i,\na_{E_k}E_j) \nonumber \\
&&=E_kE_jf_i-E_k\<f,D_{E_j}E_i\>+E_k(\Pi(E_i,E_j)\<f,\vec{\mathbf{n}}\>)-\na f(\na_{E_k}E_i,E_j) \nonumber \\
&&-\na f(E_i,\na_{E_k}E_j)\label{eq3.15} \\
&&=D^2f_i(E_j,E_k)+D_{E_k}E_j(f_i)-E_k\<f,D_{E_j}E_i\>+E_k(\Pi(E_i,E_j)\<f,\vec{\mathbf{n}}\>) \nonumber \\
&&-\na f(D_{E_k}E_i,E_j)+\Pi(E_i,E_k)\na f(\vec{\mathbf{n}},E_j)-\na f(E_i,D_{E_k}E_j)\label{eq3.16}
\eeq
where by $\na_{E_i}E_j=D_{E_i}E_j-\Pi(E_i,E_j)\vec{\mathbf{n}},$ we have
\begin{eqnarray}
&&\na f(E_i,E_j)=\<\na_{E_j}f,E_i\>=E_jf_i-\<f,\na_{E_j}E_i\> \nonumber\\
&&=E_jf_i-\<f,D_{E_j}E_i\>+\Pi(E_i,E_j)\<f,\vec{\mathbf{n}}\>.\label{eq3.17}
\end{eqnarray}
Now we compute $\na(\sym\na f)(E_i,E_j,E_k):$
\beq
&&\na(\sym\na f)(E_i,E_j,E_k)=\na_{E_k}(\sym\na f)(E_i,E_j)=E_k(\sym\na f(E_i,E_j)) \nonumber \\
&&-\sym\na f(\na_{E_k}E_i,E_j)-\sym\na f(E_i,\na_{E_k}E_j) \label{eq3.18} \\
&&=\frac{1}{2}E_kE_jf_i+\frac{1}{2}E_kE_if_j-\sym\na f(\na_{E_k}E_i,E_j)-\sym\na f(E_i,\na_{E_k}E_j) \nonumber \\
&&-\frac{1}{2}E_k(2\<f,D_{E_j}E_i\>+\<f,[E_i,E_j]\>-2\Pi(E_i,E_j)\<f,\vec{\mathbf{n}}\>)  \nonumber\\
&&=\frac{1}{2}D^2f_i(E_j,E_k)+\frac{1}{2}D_{E_k}E_j(f_i)+\frac{1}{2}D^2f_j(E_i,E_k)+\frac{1}{2}D_{E_k}E_i(f_j)                                                                                                   \nonumber \\
&&-\sym\na f(D_{E_k}E_i,E_j)+\frac{1}{2}\Pi(E_k,E_i)\<\na_{E_j}f,\vec{\mathbf{n}}\>-\sym\na f(D_{E_k}E_j,E_i) \nonumber \\
&&+\frac{1}{2}\Pi(E_k,E_j)\<\na_{E_i}f,\vec{\mathbf{n}}\>-\frac{1}{2}E_k(2\<f,D_{E_j}E_i\>+\<f,[E_i,E_j]\>-2\Pi(E_i,E_j)\<f,\vec{\mathbf{n}}\>).\label{eq3.19}
\eeq
By replacing the order of index in \eqref{eq3.19}, we can obtain $\na(\sym\na f)(E_k,E_i,E_j)$ and \\
$\na(\sym\na f)(E_j,E_k,E_i).$
Then, by \eqref{eq3.19}, we have
\beq
&&\na(\sym\na f)(E_i,E_j,E_k)+\na(\sym\na f)(E_k,E_i,E_j)-\na(\sym\na f)(E_j,E_k,E_i) \nonumber \\
&&=D^2f_i(E_j,E_k)+D_{E_k}E_j(f_i)+\frac{1}{2}[E_k,E_i](f_j)+\frac{1}{2}[E_j,E_k](f_i) \nonumber \\
&&-\sym\na f([E_k,E_i],E_j)-2\sym\na f(D_{E_k}E_j,E_i)-\sym\na f([E_j,E_k],E_i) \nonumber \\
&&-\sym\na f([E_j,E_i],E_k)+\Pi(E_k,E_j)\<\na_{E_i}f,\vec{\mathbf{n}}\> \nonumber \\
&&-\frac{1}{2}E_k(2\<f,D_{E_j}E_i\>+\<f,[E_i,E_j]\>-2\Pi(E_i,E_j)\<f,\vec{\mathbf{n}}\>) \nonumber \\
&&-\frac{1}{2}E_j(2\<f,D_{E_i}E_k\>+\<f,[E_k,E_i]\>-2\Pi(E_i,E_k)\<f,\vec{\mathbf{n}}\>) \nonumber \\
&&+\frac{1}{2}E_i(2\<f,D_{E_k}E_j\>+\<f,[E_j,E_k]\>-2\Pi(E_j,E_k)\<f,\vec{\mathbf{n}}\>), \label{eq3.20}
\eeq
where we make use of the symmetry of the Hessian of functions on $S.$
Submit \eqref{eq3.16} into \eqref{eq3.20}, we obtain
\begin{eqnarray*}
&&\na(\sym\na f)(E_i,E_j,E_k)+\na(\sym\na f)(E_k,E_i,E_j)-\na(\sym\na f)(E_j,E_k,E_i)  \\
&&=\na^2f(E_i,E_j,E_k)+\na f(D_{E_k}E_i,E_j)-\na f(D_{E_k}E_j,E_i)-\Pi(E_i,E_k)\na f(\vec{\mathbf{n}},E_j) \\
&&+\frac{1}{2}[E_j,E_k](f_i)+\frac{1}{2}[E_k,E_i](f_j)-\sym\na f([E_k,E_i],E_j)-\sym\na f([E_j,E_k],E_i) \\
&&+\Pi(E_k,E_j)\<\na_{E_i}f,\vec{\mathbf{n}}\>-\sym\na f([E_j,E_i],E_k)-\frac{1}{2}E_k\<f,[E_i,E_j]\> \\
&&-\frac{1}{2}E_j(2\<f,D_{E_i}E_k\>+\<f,[E_k,E_i]\>-2\Pi(E_i,E_k)\<f,\vec{\mathbf{n}}\>) \nonumber \\
&&+\frac{1}{2}E_i(2\<f,D_{E_k}E_j\>+\<f,[E_j,E_k]\>-2\Pi(E_j,E_k)\<f,\vec{\mathbf{n}}\>)
\end{eqnarray*}
Hence, by
\begin{eqnarray*}
&&\na f(D_{E_k}E_i,E_j)=\<\na_{E_j}f,D_{E_k}E_i\>=E_j\<f,D_{E_k}E_i\>-\<f,\na_{E_j}(D_{E_k}E_i)\>\quad\mbox{and} \\
&&\Pi(E_k,E_j)\<\na_{E_i}f,\vec{\mathbf{n}}\>=\Pi(E_k,E_j)(E_j\<f,\vec{\mathbf{n}}\>-\Pi(f^{\top},E_j)),
\end{eqnarray*}
we arrive at \eqref{eq3.14}.  \hfill$\Box$

Now we start to cope with the proof of Proposition \ref{pro3.2}.

{\bf Proof of Proposition \ref{pro3.2}}\,\,\, As we have shown before, $\sym\na y=\tilde{U}$ is equivalent with the system \eqref{eq3.8} and \eqref{eq3.9} satisfied by new variable $v,u,$ which are defined by \eqref{eq3.6} and \eqref{eq3.7}, respectively. When $S$ is elliptic, \eqref{eq3.6} is actually \eqref{eq3.11}.

Now we consider the elliptic equation \eqref{eq3.11} satisfied by $v$ as follows:
\begin{equation*}
\triangle_{\hat{g}}v=\frac{1}{\kappa}P(\tilde{U})-vtr_g\Pi+\frac{1}{2\kappa}\tilde{X}(v),
\end{equation*}
where $\hat{g}=\Pi.$
Since $S=y_0(\Om)$ and $y_0$ is an smooth isometry between $S$ and $\Om,$ \eqref{eq3.11} is also an elliptic equation in $\Om\subset\mathbb{R}^2.$

By the assumptions in Proposition \ref{pro3.2} for elliptic surface $S,$ we have that $\Pi,tr_g\Pi,\kappa\in\mathcal{C}^{2,\alpha},$ from which we see that the coefficients of the vector field
$\tilde{X}$ belong to $\mathcal{C}^{1,\alpha}.$ And since $\tilde{U}=\sym((\na\phi)^T\na \psi)\in\mathcal{C}^{2,\alpha},$ we obtain that $P(\tilde{U})\in\mathcal{C}^{0,\alpha}.$ They remain the regularity when they are projected on $\Om$ through the parametrization $y_0.$ Thus, we follow the methods in \cite{LMP} and extend the coefficients to the larger domain $\Om_{\varepsilon}=\{x'\in\mathbb{R}^2: dist(x',\Om)<\varepsilon\},$ where $\varepsilon>0$ is a sufficiently small given positive number. The extension(refer to Lemma 6.37 in \cite{GT}) preserves the regularity of these coefficients and the norms of the extension of these coefficients in $\Om_{\varepsilon}$ can be controlled by the norms of these coefficients in $\Om.$ Especially, the extension also preserves the sign of $tr_g\Pi$ and $\kappa,$ which means that \eqref{eq3.11} is still elliptic in $\Om_{\varepsilon}.$
Therefore, we consider the Dirichlet problem of elliptic equation \eqref{eq3.11} in $\Om_{\varepsilon}$ with its coefficients extended in $\Om_{\varepsilon}$ and $v=0$ on $\pl\Om_{\varepsilon}.$

By Proposition \ref{pro3.1}, we have a solution $v\in W^{2,2}(\Om_{\varepsilon})$ to \eqref{eq3.11} in $\Om_{\varepsilon}.$  Using the Sobolev embedding theorem, we obtain that $v\in \mathcal{C}^{2,\alpha}(\bar{\Om}_{\varepsilon}).$ Thus, according to the interior estimates in Schauder theory for elliptic equations, we have
\be\label{eq3.21}
\|v\|_{\mathcal{C}^{2,\alpha}(\bar{\Om})}\leqslant C\|P(\tilde{U})\|_{\mathcal{C}^{0,\alpha}(\bar{\Om})}+C\|v\|_{L^{\infty}(\Om)}\leqslant C\|\tilde{U}\|_{\mathcal{C}^{2,\alpha}(\bar{\Om})}+C\|v\|_{L^{\infty}(\Om)}.
\ee
Then we apply Theorem 8.15 in \cite{GT} to \eqref{eq3.11} in $\Om_{\varepsilon}$ and attain that
\be\label{eq3.22}
\sup_{\Om_{\varepsilon}}|v|\leqslant C(\|v\|_{L^2(\Om_{\varepsilon})}+\|\tilde{U}\|_{\mathcal{C}^{2,\alpha}(\bar{\Om})}).
\ee
By Proposition \ref{pro3.1}, we have
\be\label{eq3.23}
\|v\|_{L^2(\Om_{\varepsilon})}\leqslant C\|\tilde{U}\|_{L^{2}(\Om_{\varepsilon}))}\leqslant C\|\tilde{U}\|_{\mathcal{C}^{2,\alpha}(\bar{\Om})}.
\ee
Combine \eqref{eq3.21}, \eqref{eq3.22} and \eqref{eq3.23} and obtain
\be\label{eq3.24}
\|v\|_{\mathcal{C}^{2,\alpha}(\bar{\Om})}\leqslant C\|\tilde{U}\|_{\mathcal{C}^{2,\alpha}(\bar{\Om})}.
\ee
Thus, we have
\be\label{eq3.24a}
\|v\|_{\mathcal{C}^{2,\alpha}(\bar{S})}\leqslant C\|\tilde{U}\|_{\mathcal{C}^{2,\alpha}(\bar{S})}.
\ee
From \eqref{eq3.24a} and \eqref{eq3.9}, we arrive at
\be\label{eq3.25}
\|u\|_{\mathcal{C}^{1,\alpha}(\bar{S})}\leqslant C\|\tilde{U}\|_{\mathcal{C}^{2,\alpha}(\bar{S})}+C\|Dv\|_{\mathcal{C}^{1,\alpha}(\bar{S})}\leqslant C\|\tilde{U}\|_{\mathcal{C}^{2,\alpha}(\bar{S})}.
\ee
Via \eqref{eq3.12}, we have
\be\label{eq3.26}
\|\na y\|_{\mathcal{C}^{1,\alpha}(\bar{S})}\leqslant C(\|\tilde{U}\|_{\mathcal{C}^{1,\alpha}(\bar{S})}+\|u\|_{\mathcal{C}^{1,\alpha}(\bar{S})}+\|v\|_{\mathcal{C}^{1,\alpha}(\bar{S})})\leqslant C\|\tilde{U}\|_{\mathcal{C}^{2,\alpha}(\bar{S})}.
\ee
By Theorem 2.1 in \cite{Yao}, we obtain
\be\label{eq3.26a}
\|y\|_{\mathcal{C}^{2,\alpha}(\bar{S})}\leqslant C\|\tilde{U}\|_{\mathcal{C}^{2,\alpha}(\bar{S})}.
\ee
We apply Lemma \ref{lem3.1a} to $y=\mathcal{T}(\tilde{U})$ and from \eqref{eq3.14} and \eqref{eq3.18}, we can see that $\na^2y(E_i,E_j,E_k)$ is presented by the combination of
the components of the vector field $y$ and its first order derivatives and the components of the 2nd order tensor field $\tilde{U}$ and its first order derivatives.
Thus, we have
\be\label{eq3.27}
\|\na^2y\|_{\mathcal{C}^{1,\alpha}(\bar{S})}\leqslant C\|\tilde{U}\|_{\mathcal{C}^{2,\alpha}(\bar{S})}+C\|y\|_{\mathcal{C}^{2,\alpha}(\bar{S})}\leqslant C\|\tilde{U}\|_{\mathcal{C}^{2,\alpha}(\bar{S})}.
\ee
Therefore, note that $\tilde{U}=\sym((\na\phi)^T\na \psi)$ and combine \eqref{eq3.26a} and \eqref{eq3.27}, we obtain \eqref{eq3.13} and complete the proof.   \hfill$\Box$
\begin{rem}\label{rem3.1}
Our proof of Proposition \ref{pro3.2} provides another alternative approach for Lemma 5.1 in \cite{LMP}. In the proof there, as $\tilde{U}=\sym((\na\phi)^T\na \psi)\in\mathcal{C}^{1,\alpha},$ we use the Schauder theory for elliptic equations in divergence form \cite{GT} Theorem 8.32 and Theorem 8.15 to obtain the $\mathcal{C}^{1,\alpha}$ estimates of $v.$ Then we obtain the $\mathcal{C}^{0,\alpha}$ estimates of $u$ by \eqref{eq3.9} and thus, the $\mathcal{C}^{1,\alpha}$ estimates of $y.$ At last, by \eqref{eq3.14} holding in distribution sense, we'll recover the $\mathcal{C}^{2,\alpha}$ estimates of $y=\mathcal{T}(\sym((\na\phi)^T\na \psi)).$ From \eqref{eq3.15} and \eqref{eq3.19}, we infer that \eqref{eq3.14} will hold in the distribution sense.
\end{rem}

Now, by \eqref{eq3.13} and Banach fixed point theorem as in \cite{LMP}, we construct the required isometry on $S.$

\begin{thm}\label{thm3.2}
Let $S$ be elliptic with the $\mathcal{C}^{5,\alpha}$ regularity up to the boundary and $\mathcal{C}^{4,\alpha}$ boundary $\pl S,$ for some $0<\alpha<1.$ Given $V\in\mathcal{V}\cap\mathcal{C}^{3,\alpha}(\bar{S}),$ there exists a sequence $w_h:\bar{S}\rightarrow\mathbb{R}^3,$ equibounded in $\mathcal{C}^{3,\alpha}(\bar{S}),$ and such that for all small $h>0$ the map $u_h=id+hV+h^2w_h$ is an (exact) isometry.
\end{thm}

Now we give theorem about the density property for infinitesimal isometry on $S.$
\begin{thm}\label{thm3.3}
Let $S$ be elliptic with the $\mathcal{C}^{m+2,\alpha}$ regularity up to the boundary and $\mathcal{C}^{m+1,\alpha}$ boundary $\pl S,$ for some $0<\alpha<1$ and an integer $m>0.$ Then for every $V\in\mathcal{V},$ there exists a sequence $V_n\in\mathcal{V}\cap\mathcal{C}^{m,\alpha}(\bar{S},\mathbb{R}^3)$ such that
\begin{equation*}
\lim_{n\rightarrow\infty}\|V_n-V\|_{W^{2,2}}(S)=0.
\end{equation*}
\end{thm}

\subsection{The results on developable surfaces}
\hskip\parindent
For the case that $S=y_0(\Om)$ is a developable surface, we obtain the results following the approaches used in \cite{HLP}.
\begin{thm}\label{thm3.4}
Assume that $S$ is developable of class $\mathcal{C}^{2,1}$ and satisfy \eqref{eq1.5}, and let $\alpha\in(0,1).$ Then there exists a constant $C$ such that the following is true. For every symmetric bilinear form $B\in\mathcal{C}^{1,1}(S,T_{sym}^2(S))$ there exists a solution $w=w^{\top}+\<w,\vec{\mathbf{n}}\>\vec{\mathbf{n}}$ with
$w^{\top}\in\mathcal{C}^{0,\alpha}$ and $\<w,\vec{\mathbf{n}}\>\in L^{\infty}$ of:
\begin{equation*}
\sym\na w=\sym Dw^{\top}+\<w,\vec{\mathbf{n}}\>\Pi=B
\end{equation*}
satisfying the bounds:
\begin{equation*}
\|w^{\top}\|_{\mathcal{C}^{0,\alpha}}+\|\<w,\vec{\mathbf{n}}\>\|_{L^{\infty}}\leqslant C\|B\|_{\mathcal{C}^{1,1}}.
\end{equation*}
If, in addition, $S\in\mathcal{C}^{k+2,1}$ and $B\in\mathcal{C}^{k+1,1}$ for some $k\geqslant1,$ then
\begin{equation*}
\|w^{\top}\|_{\mathcal{C}^{k,1}}+\|\<w,\vec{\mathbf{n}}\>\|_{\mathcal{C}^{k-1,1}}\leqslant C\|B\|_{\mathcal{C}^{k+1,1}}.
\end{equation*}
\end{thm}

Now we give the definition of mth order infinitesimal isometry.
\begin{dfn}\label{def3.1}
An one parameter family $\{u_{\epsilon}\}_{\epsilon>0}\subset\mathcal{C}^{0,1}(\bar{S},\mathbb{R}^3)$ is said to be a (generalized) mth order infinitesimal isometry if the change of metric induced by $u_{\epsilon>0}$ is of order $\epsilon^{m+1},$ that is
\begin{equation*}
\|(\na u_{\epsilon})^T\na u_{\epsilon}-g\|_{L^{\infty}(S)}=O(\epsilon^{m+1})\quad\mbox{as}\quad \epsilon\rightarrow0,
\end{equation*}
where $g$ is the induced metric on $S\subset\mathbb{R}^3.$
\end{dfn}
Then, we give the results of matching property and density of infinitesimal isometry. Because we need the mth order infinitesimal isometry with higher regularity than that in \cite{HLP} in the construction of the recovery sequence, we apply the iteration procedure(that is Lemma \ref{lem3.1} below) proved in Theorem 5.2 in \cite{HLP} to produce the mth order infinitesimal isometry we need.
\begin{lem}\label{lem3.1}
Let $S$ be of class $\mathcal{C}^{k+2,1},$ where $k\in\mathbb{N},$ and let $u_{\epsilon}$ be an $(i-1)$th order isometry of regularity $\mathcal{C}^{k+1,1}$ of the form:
\begin{equation*}
u_{\epsilon}=id+\sum_{j=1}^{i-1}\epsilon^jw_j,\quad w_j\in\mathcal{C}^{k+1,1}.
\end{equation*}
Then there exists $w_i\in\mathcal{C}^{k-1,1}(S,\mathbb{R}^3)$ so that $\phi_{\epsilon}=u_{\epsilon}+\epsilon^iw_i$ is an ith order infinitesimal isometry, and
\begin{equation*}
\|w_i\|_{\mathcal{C}^{k-1,1}}\leqslant C\sum_{j=1}^{i-1}\|w_j\|_{\mathcal{C}^{k+1,1}}\|w_{i-j}\|_{\mathcal{C}^{k+1,1}}.
\end{equation*}
\end{lem}
Applying this lemma, we have the following theorem.
\begin{thm}\label{thm3.5}
Let $S$ be developable and satisfy $\bf(H_d).$ Given $V\in\mathcal{V}\cap\mathcal{C}^{2m+1,1}(\bar{S}),$ there exists a sequence $w_{\epsilon}:\bar{S}\rightarrow\mathbb{R}^3,$ equibounded in $\mathcal{C}^{3,1}(S),$ and such that for all small $\epsilon>0$ the map $u_{\epsilon}=id+\epsilon V+\epsilon^2w_{\epsilon}$ is a (generalized) mth order infinitesimal isometry of class $\mathcal{C}^{3,1}.$
\end{thm}
\begin{thm}\label{thm3.6}
Let $S$ be developable and satisfy $\bf(H_d)$ with the $\mathcal{C}^{k+1,1}$ regularity up to the boundary. Then for every $V\in\mathcal{V},$ there exists a sequence $V_n\in\mathcal{V}\cap\mathcal{C}^{k,1}(\bar{S},\mathbb{R}^3)$ such that
\begin{equation*}
\lim_{n\rightarrow\infty}\|V_n-V\|_{W^{2,2}}(S)=0.
\end{equation*}
\end{thm}

\subsection{The results on hyperbolic surfaces}
\hskip\parindent
For the case that $S=y_0(\Om)$ is a hyperbolic surface, we have the results following \cite{Yao}.
\begin{thm}\label{thm3.7}
Assume that $S$ is hyperbolic and satisfy $\bf(H_h)$ with $\mathcal{C}^{2,1}$ regularity. For every symmetric bilinear form $\tilde{U}\in\mathcal{C}^{1,1}(S,T_{sym}^2(S))$ there exists a solution $y=y^{\top}+\<y,\vec{\mathbf{n}}\>\vec{\mathbf{n}}\in\mathcal{C}^{0,1}(S,\mathbb{R}^3)$ to \eqref{eq3.5} satisfying the bounds
\begin{equation*}
\|y^{\top}\|_{\mathcal{C}^{1,1}}+\|\<y,\vec{\mathbf{n}}\>\|_{\mathcal{C}^{0,1}}\leqslant C\|\tilde{U}\|_{\mathcal{C}^{1,1}}.
\end{equation*}
If, in addition, $S\in\mathcal{C}^{m+2,1}$ and $B\in\mathcal{C}^{m+1,1}$ for some $m\geqslant1,$ then
\begin{equation*}
\|y^{\top}\|_{\mathcal{C}^{m+1,1}}+\|\<y,\vec{\mathbf{n}}\>\|_{\mathcal{C}^{m,1}}\leqslant C\|\tilde{U}\|_{\mathcal{C}^{m+1,1}}.
\end{equation*}
\end{thm}

Now we turn to the matching property and density results to infinitesimal isometry on hyperbolic surfaces. Similarly, by applying Lemma \ref{lem3.1}, we have the following theorem to produce the mth order infinitesimal isometry with higher regularity than that in \cite{Yao}.
\begin{thm}\label{thm3.8}
Let $S$ be hyperbolic and satisfy $\bf(H_h).$ Given $V\in\mathcal{V}\cap\mathcal{C}^{2m+1,1}(\bar{S}),$ there exists a sequence $w_{\epsilon}:\bar{S}\rightarrow\mathbb{R}^3,$ equibounded in $\mathcal{C}^{3,1}(S),$ and such that for all small $\epsilon>0,$ $u_{\epsilon}=id+\epsilon V+\epsilon^2w_{\epsilon}$ is a (generalized) mth order infinitesimal isometry of class $\mathcal{C}^{3,1}.$
\end{thm}
\begin{thm}\label{thm3.9}
Let $S$ be hyperbolic and satisfy $\bf(H_h)$ with $\mathcal{C}^{m+2,1}$ regularity for some integer $m\geqslant0.$ Then for every $V\in\mathcal{V},$ there exists a sequence $V_n\in\mathcal{V}\cap\mathcal{C}^{m,1}(\bar{S},\mathbb{R}^3)$ such that
\begin{equation*}
\lim_{n\rightarrow\infty}\|V_n-V\|_{W^{2,2}}(S)=0.
\end{equation*}
\end{thm}

Therefore, we use the $1-1$ corresponding between the infinitesimal isometry $\tilde{V}$ on $S$ and $V$ in \eqref{eq2.4b} to pull back the constructions on surface $S$ above to the mid-plate $(\Om,g).$ Then we are prepared to construct the recovery sequence.

\section{The upper bound of the scaled elastic energy}
\setcounter{equation}{0}
\hskip\parindent
In this Section, we construct the recovery sequence and obtain the upper bound of the elastic energy with incompatible prestrain, which finish the computation of the $\Ga-$ limit.

{\bf Proof of Theorem \ref{thm1.2}}\,\,\, By density results and the continuity of $I_{\beta}$ with respect to the strong topology of $W^{2,2},$ we can assume that $V\in\mathcal{V}\cap\mathcal{C}^{3,\alpha}(\bar{S},\mathbb{R}^3).$ In general case, the results will then follow from a diagonal arguments. We denote $S=y_0(\Om).$

Step 1.\,\,\, Let $\epsilon=\frac{\sqrt{e^h}}{h}.$ We recall that $\epsilon\rightarrow0$ as $h\rightarrow0,$ by assumption \eqref{eq1.1}. Therefore, by Theorem \ref{thm3.2}, there
exists a sequence $w_{\epsilon}:\bar{S}\rightarrow\mathbb{R}^3,$ equibounded in $\mathcal{C}^{3,\alpha}(\bar{S}),$ such that for all small $h>0$ the map
\begin{equation*}
u_{\epsilon}=id+\epsilon V+\epsilon^2w_{\epsilon}
\end{equation*}
is an exact isometry on the surface $S=y_0(\Om).$ Thus, by the 1-1 corresponding of the vector fields between $\Om$ and $S,$ we have
\be\label{eq4.1}
u_{\epsilon}(y_0(x'))=y_0(x')+\epsilon V(y_0(x'))+\epsilon^2w_{\epsilon}(y_0(x'))
\ee
is alao an exact isometry on the two-dimensional Riemannian manifold $(\Om,g=G_{2\times2}).$

For every $x'\in\Om,$ let $\vec{\mathbf{n}}_{\epsilon}(x')$ denote the unit vector normal to $u_{\epsilon}(S)$ at the point $u_{\epsilon}(y_0(x')).$ By the regularity of $u_{\epsilon}$ we have that $\vec{\mathbf{n}}_{\epsilon}\in\mathcal{C}^{2,\alpha}(\bar{S},\mathbb{R}^3),$ while by \eqref{eq4.1} and $V\in\mathcal{V}$ we obtain the expansion
\beq\label{eq4.2}
&&\vec{\mathbf{n}}_{\epsilon}(y_0(x'))=\nabla_{E_1}u_{\epsilon}\times\nabla_{E_2}u_{\epsilon}
=(E_1+\epsilon\nabla_{E_1}V+\epsilon^2\nabla_{E_1}w_{\epsilon})\times(E_2+\epsilon\nabla_{E_2}V+\epsilon^2\nabla_{E_2}w_{\epsilon}) \nonumber\\
&&=E_1\times E_2+\epsilon(E_1\times\nabla_{E_2}V+\nabla_{E_1}V\times E_2) \nonumber\\
&&+\epsilon^2(E_1\times\nabla_{E_2}w_{\epsilon}+\nabla_{E_1}V\times\nabla_{E_2}V+\nabla_{E_1}w_{\epsilon}\times E_2)  \nonumber\\
&&+\epsilon^3(\nabla_{E_1}V\times\nabla_{E_2}w_{\epsilon}+\nabla_{E_1}w_{\epsilon}\times\nabla_{E_2}V)+\epsilon^4\nabla_{E_1}w_{\epsilon}\times\nabla_{E_2}w_{\epsilon} \nonumber\\
&&=\vec{\mathbf{n}}(y_0(x'))+\epsilon A\vec{\mathbf{n}}(y_0(x'))+O(\epsilon^2).
\eeq
Indeed, one can take the following form, which simplifies the computations
\begin{equation*}
\vec{\mathbf{n}}_{\epsilon}(y_0(x'))=\frac{\pl_1u_{\epsilon}(y_0(x'))\times\pl_2u_{\epsilon}(y_0(x'))}{|\pl_1u_{\epsilon}(y_0(x'))\times\pl_2u_{\epsilon}(y_0(x'))|}
=\nabla_{E_1}u_{\epsilon}\times\nabla_{E_2}u_{\epsilon},
\end{equation*}
where $\{E_1,E_2\}$ is an orthonormal vector field on $S=y_0(\Om)$ with the same orientation with $\{\pl_1y_0,\pl_2y_0\}.$ This means that
\begin{equation*}
\vec{\mathbf{n}}=\frac{\pl_1y_0(x')\times\pl_2y_0(x')}{|\pl_1y_0(x')\times\pl_2y_0(x')|}
=E_1\times E_2.
\end{equation*}

Now, we define the vector field $\vec{b}_{\epsilon}(y_0(x'))$ on the surface $S$ satisfying
\be\label{eq4.3}
Q_{\epsilon}=[\pl_1u_{\epsilon},\pl_2u_{\epsilon},\vec{b}_{\epsilon}],\quad(Q_{\epsilon})^TQ_{\epsilon}=G(x')\quad\mbox{and}\quad \det Q_{\epsilon} >0.
\ee
By the methods used in \cite{BLS}, we arrive at
\be\label{eq4.4}
\vec{b}_{\epsilon}(x')=-\frac{1}{G^{33}}(G^{13}\pl_1u_{\epsilon}+G^{23}\pl_2u_{\epsilon})+\frac{1}{\sqrt{G^{33}}}\vec{\mathbf{n}}_{\epsilon}(y_0(x')).
\ee
Moreover, by \eqref{eq4.2} and \eqref{eq4.1}, we have
\beq\label{eq4.5}
&&\vec{b}_{\epsilon}(x')=-\frac{1}{G^{33}}(G^{13}\pl_1u_{\epsilon}+G^{23}\pl_2u_{\epsilon})+\frac{1}{\sqrt{G^{33}}}\vec{\mathbf{n}}_{\epsilon}(y_0(x')) \nonumber \\
&&=\vec{b}_0(x')+\epsilon[-\frac{1}{G^{33}}(G^{13}\pl_1V+G^{23}\pl_2V)+\frac{1}{\sqrt{G^{33}}}A\vec{\mathbf{n}}(y_0(x'))]+O(\epsilon^2) \nonumber \\
&&=\vec{b}_0(x')+\epsilon A\vec{b}_0(x')+O(\epsilon^2),
\eeq
where we used the definition of the space of infinitesimal isometry $\mathcal{V}.$

Now we define the vector fields $d^h(y_0(x'))$ and $d_{\epsilon}(y_0(x')).$
First we define $d^h(y_0(x'))\in W^{1,\infty}(S,\mathbb{R}^3)$ such that
\beq
&&\lim_{h\rightarrow0}\sqrt{h}\|d^h\|_{W^{1,\infty}}=0, \label{eq4.6}\\
&&\lim_{h\rightarrow0}d^h=(Q_0^T)^{-1}c(x',(\nabla y_0)^T\nabla\vec{p}+(\nabla V)^T\nabla\vec{b}_0))\quad\mbox{in}\quad L^{\infty}(\Om), \label{eq4.6a}
\eeq
where $c(x',F_{2\times2})$ is the unique minimizer in \eqref{eq2.4c}. Note that $c(x',\cdot)$ is a linear function of $F_{2\times2}$ and it only depends on its symmetric part
$(\sym F_{2\times2}).$ Moreover, by the regularity of the vector field $V,$ the vector field $d(x')=(Q_0^T)^{-1}c(x',(\nabla y_0)^T\nabla\vec{p}+(\nabla V)^T\nabla\vec{b}_0))$ belongs to $L^{\infty}.$
We then define another vector field $d_{\epsilon}(y_0(x'))\in\mathcal{C}^{1,\alpha}(\bar{S},\mathbb{R}^3)$ such that
\be\label{eq4.7}
(\na u_{\epsilon})^Td_{\epsilon}=-(\na\vec{b}_{\epsilon})^T\cdot\vec{b}_{\epsilon},\quad\mbox{and}\quad\<\vec{b}_{\epsilon},d_{\epsilon}\>=0.
\ee
By \eqref{eq4.7}, we have
\begin{equation*}
d_{\epsilon}=-(Q_{\epsilon}^T)^{-1}(\<\pl_1\vec{b}_{\epsilon},\vec{b}_{\epsilon}\>,\<\pl_2\vec{b}_{\epsilon},\vec{b}_{\epsilon}\>,0)^T,
\end{equation*}
which leads to the fact that $d_{\epsilon}\rightarrow\vec{d}_0$ in the strong topology of $L^{\infty}(\Om).$

Here we introduce the recovery sequence $u^h$ as required by the statement of the theorem. Note that the following suggestion for $u^h$ is in accordance with the one used in
\cite{FJ2M} in the framework of the purely nonlinear bending theory for shells, corresponding to the scaling regime $\beta=2.$ Consider the sequence of deformations $u^h\in W^{1,2}(\Om^h,\mathbb{R}^3)$ defined by
\be\label{eq4.8}
u^h(x',x_3)=u_{\epsilon}(y_0(x'))+x_3\vec{b}_{\epsilon}(y_0(x'))+\frac{x_3^2}{2}d_{\epsilon}(y_0(x'))+\frac{x_3^2}{2}\epsilon d^h(y_0(x'))\quad\mbox{on}\quad \Om^h
\ee
By change of variable, we have
\beq\label{eq4.9}
&&y^h(x',x_3)=u^h(x',hx_3)=u_{\epsilon}(y_0(x'))+hx_3\vec{b}_{\epsilon}(y_0(x'))+\frac{h^2x_3^2}{2}d_{\epsilon}(y_0(x')) \nonumber \\
&&+\frac{h^2x_3^2}{2}\epsilon d^h(y_0(x'))\quad\mbox{on}\quad \Om^1.
\eeq
Therefore, properties $(i),(ii),(iii)$ mentioned in Theorem \ref{thm1.2} now easily follow from the uniform bound on $w_{\epsilon}$ and \eqref{eq4.6}.

Step 2.\,\,\,To prove \eqref{eq1.5a}, it's convenient to perform a change of variable in the energy $E^h(u^h)$ and express it in terms of the scaled deformation $y^h.$ By straight calculation, we have
\be\label{eq4.10}
\frac{1}{e^h}E^h(u^h)=\frac{1}{e^h}\int_{-\frac{1}{2}}^{\frac{1}{2}}\int_{\Om}W(\na_h y^h(x',x_3)A^{-1})dx'dx_3,
\ee
where $\na_h y^h(x',x_3)=\na u^h(x',hx_3).$ We also have
\beq
&&\frac{1}{h}\pl_3y^h(x',x_3)=\vec{b}_{\epsilon}(y_0(x'))+hx_3d_{\epsilon}(y_0(x'))+hx_3\epsilon d^h(y_0(x')) \label{eq4.11} \\
&&\pl_iy^h(x',x_3)=\pl_i(u_{\epsilon}\circ y_0)(x')+hx_3\pl_i(\vec{b}_{\epsilon}\circ y_0)(x')+\frac{h^2x_3^2}{2}\pl_i(d_{\epsilon}\circ y_0)(x') \nonumber \\
&&+\frac{h^2x_3^2}{2}\epsilon \pl_i(d^h\circ y_0)(x')\quad i=1,2\label{eq4.12}
\eeq

From \eqref{eq4.1}, \eqref{eq4.5} and \eqref{eq4.6}, it follows that
\be\label{eq4.13}
\|\na_hy^h-Q_0\|_{L^{\infty}}\rightarrow0,\quad\mbox{as} h\rightarrow0.
\ee
It now follows by polar decomposition theorem (for $h$ sufficiently small) that there exists a proper rotation $R_0\in SO(3)$ and a well defined square root of $A^{-1}(\na_hy^h)^T\na_hy^hA^{-1}$ such that
\begin{equation*}
\na_hy^hA^{-1}=R_0\sqrt{A^{-1}(\na_hy^h)^T\na_hy^hA^{-1}}.
\end{equation*}
By frame difference of $W,$ we deduce that
\begin{equation*}
W(\na_hy^hA^{-1})=W(\sqrt{A^{-1}(\na_hy^h)^T\na_hy^hA^{-1}})=W(id+\frac{1}{2}A^{-1}K^hA^{-1}+O(|K^h|^2)),
\end{equation*}
where the last equality follows by Taylor expansion, with $K^h$ given by
\begin{equation*}
K^h=(\na_hy^h)^T\na_hy^h-G.
\end{equation*}
As $\|K^h\|_{L^{\infty}}$ is infinitesimal as $h\rightarrow0,$ we use the formula
\begin{equation*}
W(id+K)=\frac{1}{2}D^2W(id)(K,K)+\int_0^1(1-s)[D^2W(id+sK)-D^2W(id)](K,K)ds
\end{equation*}
and obtain
\begin{equation}\label{eq4.14}
\frac{1}{e^h}W(\na_hy^hA^{-1})=\frac{1}{2}\mathcal{Q}_3(\frac{1}{2\sqrt{e^h}}A^{-1}K^hA^{-1}+\frac{1}{\sqrt{e^h}}O(|K^h|^2))+\frac{1}{e^h}o(|K^h|^2).
\end{equation}

Using \eqref{eq4.11} and \eqref{eq4.12}, we now calculate $K^h.$ We first consider $(K^h)_{2\times2}:$
\be\label{eq4.15}
K^h(e_i,e_j)=\<\pl_iy^h,\pl_jy^h\>-G_{ij}=2hx_3\sym((\na u_{\epsilon})^T\na\vec{b}_{\epsilon})(e_i,e_j)+o(\sqrt{e^h}),\quad i,j=1,2,
\ee
where we have used the fact that $u_{\epsilon}$ is an isometry.
To compute $\sym((\na u_{\epsilon})^T\na\vec{b}_{\epsilon})(e_i,e_j),$ we need to calculate
\beq\label{eq4.16}
&&\<\pl_iu_{\epsilon},\pl_j\vec{b}_{\epsilon}\>=\<\pl_iy_0+\epsilon \pl_i(V\circ y_0)+\epsilon^2\pl_i(w_{\epsilon}\circ y_0),\pl_j\vec{b}_0+\epsilon \pl_j(A\vec{b}_0)+O(\epsilon^2)\> \nonumber \\
&&=\<\pl_iy_0,\pl_j\vec{b}_0\>+\epsilon(\<\pl_iy_0,\pl_j\vec{p}\>+\<\pl_i(V\circ y_0),\pl_j\vec{b}_0\>)+O(\epsilon^2),
\eeq
where we use $A\vec{b}_0=\vec{p}$ which can be obtained by the definition of the vector field $\vec{p}.$
Thus, by \eqref{eq1.3}, we arrive at
\beq\label{eq4.17}
&&K^h(e_i,e_j)=\<\pl_iy^h,\pl_jy^h\>-G_{ij}=2\sqrt{e^h}x_3\sym((\na y_0)^T\na\vec{p}+(\na (V\circ y_0))^T\na\vec{b}_0)(e_i,e_j) \nonumber \\
&&+o(\sqrt{e^h}),\quad i,j=1,2,
\eeq
Now we calculate $(K^h)_{33}$ and by \eqref{eq4.11}, obtain
\beq\label{eq4.18}
&&K^h(e_3,e_3)=\frac{1}{h^2}\<\pl_3y^h,\pl_3y^h\>-G_{33}=2hx_3\epsilon\<\vec{b}_{\epsilon},d^h\>+o(\sqrt{e^h}) \nonumber \\
&&=2\sqrt{e^h}x_3\<\vec{b}_{\epsilon},d^h\>+o(\sqrt{e^h}).
\eeq
The remaining coefficients of the symmetric matrix $K^h$ are
\be\label{eq4.19}
K^h(e_i,e_3)=\frac{1}{h}\<\pl_iy^h,\pl_3y^h\>-G_{i3}
=\sqrt{e^h}x_3\<\pl_iu_{\epsilon},d^h\>+o(\sqrt{e^h}), i=1,2,
\ee
where we make use of the definition \eqref{eq4.7} of the vector field $d_{\epsilon}.$

Step 3. From the previous calculations \eqref{eq4.17}-\eqref{eq4.19}, we finally deduce that
\beq\label{eq4.20}
&&\lim_{h\rightarrow0}\frac{1}{2\sqrt{e^h}}K^h=(x_3\sym((\na y_0)^T\na\vec{p}+(\na (V\circ y_0))^T\na\vec{b}_0))^{\ast}+x_3\<\vec{b}_0,d\>e_3\otimes e_3 \nonumber \\
&&+\frac{1}{2}x_3\<d,\pl_iy_0\>e_3\otimes e_i+\frac{1}{2}x_3\<\pl_iy_0,d\>e_i\otimes e_3  \nonumber \\
&&=(x_3\sym((\na y_0)^T\na\vec{p}+(\na (V\circ y_0))^T\na\vec{b}_0))^{\ast}+x_3\sym(Q_0^Td\otimes e_3) \nonumber \\
&&=(x_3\sym((\na y_0)^T\na\vec{p}+(\na (V\circ y_0))^T\na\vec{b}_0))^{\ast} \nonumber \\
&&+x_3\sym(c(x',(\nabla y_0)^T\nabla\vec{p}+(\nabla V)^T\nabla\vec{b}_0)\otimes e_3) \quad\mbox{in}\quad L^{\infty}(\Om),
\eeq
where the symbol $(F)^{\ast}$ stands for the matrix such that $((F)^{\ast})_{2\times2}=F_{2\times2}$ and other elements of $(F)^{\ast}$ are all zero.
Using \eqref{eq4.10}, \eqref{eq4.14}, \eqref{eq4.20} and the dominated convergence theorem, we obtain
\begin{eqnarray*}
&&\lim_{h\rightarrow0}\frac{1}{e^h}E^h(u^h)=\lim_{h\rightarrow0}\frac{1}{e^h}\int_{-\frac{1}{2}}^{\frac{1}{2}}\int_{\Om}W(\na_h y^h(x',x_3)A^{-1})dx'dx_3 \\
&&=\lim_{h\rightarrow0}\frac{1}{2}\int_{-\frac{1}{2}}^{\frac{1}{2}}\int_{\Om}\mathcal{Q}_3(\frac{1}{2\sqrt{e^h}}A^{-1}K^hA^{-1}+\frac{1}{\sqrt{e^h}}O(|K^h|^2))dx'dx_3 \\
&&=\frac{1}{2}\int_{-\frac{1}{2}}^{\frac{1}{2}}\int_{\Om}\mathcal{Q}_3(x_3A^{-1}((\sym((\na y_0)^T\na\vec{p}+(\na (V\circ y_0))^T\na\vec{b}_0))^{\ast} \\
&&+\sym(c(x',(\nabla y_0)^T\nabla\vec{p}+(\nabla V)^T\nabla\vec{b}_0)\otimes e_3))A^{-1})dx'dx_3 \\
&&=\frac{1}{24}\int_{\Omega}\mathcal{Q}_{2,A}(x',(\nabla y_0)^T\nabla\vec{p}+(\nabla V)^T\nabla\vec{b_0})dx'=\mathcal{I}_{\beta}(V).
\end{eqnarray*}
Thus, the proof of Theorem \ref{thm1.2} is completed.   \hfill$\Box$

Similarly to the proof of Theorem \ref{thm1.2}, we now give the proof of Theorem \ref{thm1.3}.

{\bf Proof of Theorem \ref{thm1.3}}\,\,\, We shall construct a recovery sequence for developable surfaces based on Theorem \ref{thm3.5} and \ref{thm3.6} or hyperbolic surfaces
based on Theorem \ref{thm3.8} and \ref{thm3.9}. Indeed, by the density result and the continuity of the functional $\mathcal{I}_{\beta}$ with respect to the strong topology of $W^{2,2}(\Om),$ we shall assume $V\in\mathcal{V}\cap\mathcal{C}^{2m+1,1}(\bar{S},\mathbb{R}^3)$ for both developable case and hyperbolic case. In the general case the results will then follow through a diagonal arguments. The following proof is similar to that of the elliptic case. So we only point out the difference between them and omit the detail.

Let $\epsilon$ be the same as that in the proof of Theorem \ref{thm1.2}, that is $\epsilon=\frac{\sqrt{e^h}}{h}.$ From \eqref{eq1.1}, we see that $\epsilon\rightarrow0,$ as $h\rightarrow0.$ Therefore, by Theorem \ref{thm3.5} and \ref{thm3.8}, there exists a sequence $w_{\epsilon}:\bar{S}\rightarrow\mathbb{R}^3,$ equibounded in $\mathcal{C}^{3,1}(\bar{S}),$ such that for all small $h>0:$
\be\label{eq4.21}
u_{\epsilon}=id+\epsilon V+\epsilon^2w_{\epsilon}
\ee
is a (generalized) mth order infinitesimal isometry on $S.$ Note that by \eqref{eq1.6} we have
\begin{equation*}
\frac{\epsilon^{m+1}}{\sqrt{e^h}}=\frac{(\sqrt{e^h})^m}{h^{m+1}}=\frac{o(h^{m+1})}{h^{m+1}}\rightarrow0,
\end{equation*}
here $\epsilon^{m+1}=o(\sqrt{e^h}).$ We may thus replace $O(\epsilon^{m+1})$ with $o(\sqrt{e^h}).$

Because $u_{\epsilon}$ is a (generalized) mth order infinitesimal isometry on $S,$ we have that
\be\label{eq4.22}
(\na u_{\epsilon})^T\na u_{\epsilon}=g+O(\epsilon^{m+1}).
\ee
Thus, $u_{\epsilon}:\bar{S}\rightarrow\mathbb{R}^3$ is an immersion, which means that $\{\pl_1u_{\epsilon},\pl_2u_{\epsilon}\}$ is linear independent.
We then define the unit normal vector
\begin{equation*}
\vec{\mathbf{n}}_{\epsilon}(y_0(x'))=\frac{\pl_1u_{\epsilon}(y_0(x'))\times\pl_2u_{\epsilon}(y_0(x'))}{|\pl_1u_{\epsilon}(y_0(x'))\times\pl_2u_{\epsilon}(y_0(x'))|}
=\frac{\nabla_{E_1}u_{\epsilon}\times\nabla_{E_2}u_{\epsilon}}{|\nabla_{E_1}u_{\epsilon}\times\nabla_{E_2}u_{\epsilon}|},
\end{equation*}
where $\{E_1,E_2\}$ is an orthonormal vector field on $S=y_0(\Om)$ with the same orientation with $\{\pl_1y_0,\pl_2y_0\}.$
By \eqref{eq4.22}, we conduct the similar calculations to that in Theorem \ref{thm1.2} for $\vec{\mathbf{n}}_{\epsilon}$ and conclude that
\beq
&&|\nabla_{E_i}u_{\epsilon}|^2=1+O(\epsilon^2), |\<\nabla_{E_1}u_{\epsilon},\nabla_{E_2}u_{\epsilon}\>|=O(\epsilon^2)\quad i=1,2, \nonumber \\
&&|\nabla_{E_1}u_{\epsilon}\times\nabla_{E_2}u_{\epsilon}|=1+O(\epsilon^2),\nonumber \\
&&\vec{\mathbf{n}}_{\epsilon}=\vec{\mathbf{n}}+\epsilon A\vec{\mathbf{n}}+O(\epsilon^2). \label{eq4.23}
\eeq

We define the vector field $\vec{b}_{\epsilon}(y_0(x'))$ on the surface $S$ by
\be\label{eq4.24}
Q_{\epsilon}=[\pl_1u_{\epsilon},\pl_2u_{\epsilon},\vec{b}_{\epsilon}],\quad(Q_{\epsilon})^TQ_{\epsilon}=G(x')+O(\epsilon^{m+1})\triangleq \tilde{G}\in\mathcal{C}^{2,1}(\bar{S})\quad\mbox{and}\quad \det Q_{\epsilon}>0.
\ee
It's obvious that $\tilde{G}$ is positive definite as $h$ sufficiently small.

By the methods used in \cite{BLS}, we arrive at
\be\label{eq4.25}
\vec{b}_{\epsilon}(x')=-\frac{1}{\tilde{G}^{33}}(\tilde{G}^{13}\pl_1u_{\epsilon}+\tilde{G}^{23}\pl_2u_{\epsilon})
+\frac{1}{\sqrt{\tilde{G}^{33}}}\vec{\mathbf{n}}_{\epsilon}(y_0(x'))\in\mathcal{C}^{2,1}(\bar{S}).
\ee
Then by \eqref{eq4.25} and the regularity of $u_{\epsilon}$ in Theorem \ref{thm3.5} and \ref{thm3.8}, we still have
\beq\label{eq4.26}
&&\vec{b}_{\epsilon}(x')=-\frac{1}{\tilde{G}^{33}}(\tilde{G}^{13}\pl_1u_{\epsilon}+\tilde{G}^{23}\pl_2u_{\epsilon})+\frac{1}{\sqrt{\tilde{G}^{33}}}\vec{\mathbf{n}}_{\epsilon}(y_0(x')) \nonumber \\
&&=-\frac{1}{G^{33}}(G^{13}\pl_1u_{\epsilon}+G^{23}\pl_2u_{\epsilon})+\frac{1}{\sqrt{G^{33}}}\vec{\mathbf{n}}_{\epsilon}(y_0(x'))+O(\epsilon^{m+1}) \nonumber \\
&&=\vec{b}_0(x')+\epsilon A\vec{b}_0(x')+O(\epsilon^2).
\eeq

We still employ the same $u^h$ with that in Theorem \ref{thm1.2}, as follows
\begin{equation*}
u^h(x',x_3)=u_{\epsilon}(y_0(x'))+x_3\vec{b}_{\epsilon}(y_0(x'))+\frac{x_3^2}{2}d_{\epsilon}(y_0(x'))+\frac{x_3^2}{2}\epsilon d^h(y_0(x'))\quad\mbox{on}\quad \Om^h,
\end{equation*}
where the vector fields $d_{\epsilon}\in\mathcal{C}^{1,1}(\bar{S})$ and $d^h$ are also defined the same as that in Theorem \ref{thm1.2} by \eqref{eq4.6}, \eqref{eq4.6a} and \eqref{eq4.7}.
As we obtain the (generalized) mth order infinitesimal isometry  $u_{\epsilon}$ with the $\mathcal{C}^{3,1}$ regularity, we can proceed just as that in Theorem \ref{thm1.2}
and also arrive at the same $\Ga-$ limit $\mathcal{I}_{\beta}(V)=\frac{1}{24}\int_{\Omega}\mathcal{Q}_{2,A}(x',(\nabla y_0)^T\nabla\vec{p}+(\nabla V)^T\nabla\vec{b_0})dx'.$
 \hfill$\Box$

\section{Discussion of the obtained functional}
\setcounter{equation}{0}
\hskip\parindent
Similar to \cite{L2R}, we give the comparison of our results with that obtained in \cite{LMP}. We'll see that the arguments in both energies are related via the
parametrization $y_0$ of the surface $y_0(\Om)$ in \eqref{eq1.3}.

Recall that when $S$ is a smooth $2d$ surface in $\mathbb{R}^3,$ the $\Gamma-$ limit of the scaled elastic energies $\frac{1}{e^h}E^h(u^h)$ on thin shell $S^h$ with mid-surface $S$ is
\be\label{eq5.1}
\tilde{\mathcal{I}}_{\beta,S}(V)=\frac{1}{24}\int_S\mathcal{Q}_2(x,(\nabla(\tilde{A}N)-\tilde{A}\Pi)_{tan})dx,
\ee
with $\tilde{e}=\frac{h}{\sqrt{e^h}}\lim_{h\rightarrow0}\nabla w^h\in W^{1,2}(S,\mathbb{R}^3)$ belongs to the finite strain space. Above, $\Pi$ stands for the second fundamental form of $S$ and $N$ is the unit normal vector to $S.$ Moreover, $\tilde{A}\in W^{1,2}(S,so(3))$ is associated with $\tilde{V}$ by $\pl_{\tau}\tilde{V}(x)=\tilde{A}(x)\tau,$ for any $\tau\in T_xS.$

In the present setting, denote $S=y_0(\Om)$ and observe that the $1-1$ corresponding between $\tilde{V}$ mentioned above and $V$ in \eqref{eq2.4b}, given by the change of variables $V=\tilde{V}\circ y_0.$ The skew symmetric tensor field $\tilde{A}$ on $T_xS$ is then uniquely given by
\be\label{eq5.2}
\pl_eV(x')=\tilde{A}(y_0(x'))\pl_ey_0\quad\mbox{and}\quad \tilde{A}\vec{b}_0=\vec{p}\quad \forall e\in\mathbb{R}^2,
\ee
and the extended strains above are related to that in Theorem \ref{thm2.1} by
\begin{equation*}
\langle\tilde{e}(y_0(x'))\pl_ey_0,\pl_ey_0\rangle=\langle e(x')e,e\rangle \quad \forall e\in\mathbb{R}^2.
\end{equation*}

When $\vec{b}_0=N,$ the arguments in the functional \eqref{eq5.1} coincides with $(\nabla y_0)^T\nabla\vec{p}+(\nabla V)^T\nabla\vec{b}_0.$ Indeed,
\begin{eqnarray*}
&&\langle(\pl_\tau\tilde{A})\vec{b}_0,\tau\rangle=\langle\pl_e(\tilde{A}\vec{b}_0),\pl_ey_0\rangle-\langle\tilde{A}\pl_e\vec{b}_0,\pl_ey_0\rangle
=\langle\pl_e\vec{p},\pl_ey_0\rangle+\langle\pl_e\vec{b}_0,\tilde{A}\pl_ey_0\rangle \\
&&=\langle(\nabla y_0)^T\nabla\vec{p}e,e\rangle-\langle(\nabla V)^T\nabla\vec{b}_0e,e\rangle,
\end{eqnarray*}
in view of \eqref{eq3.2}, where $\tau=\pl_ey_0\in T_{y_0(x')}S,$ for any $e\in\mathbb{R}^2.$ But in general,
\begin{equation*}
\vec{b}_0=-\frac{1}{G^{33}}(G^{13}\pl_1y_0+G^{23}\pl_2y_0)+\frac{1}{\sqrt{G^{33}}}N.
\end{equation*}
Obviously, $\vec{b}_0\neq N$ in general case on the surfsce $S.$ Thus, this means that the arguments in the functional \eqref{eq5.1} obtained in \cite{LMP} are different from that in the limit functional \eqref{eq2.5} obtained in this paper.

\end{document}